\documentclass{article}

\usepackage{amsmath}
\usepackage{amssymb}

\usepackage[2cell,v2]{xy}
\UseTwocells
\UseHalfTwocells

\newcommand{\tmem}[1]{{\em #1\/}}
\newcommand{\tmmathbf}[1]{\mathbf{#1}}
\newcommand{\tmop}[1]{{#1}}

\newcommand{\dueto}[1]{\textup{{(#1) }}}

\newcommand{\leirom}{\renewcommand{\labelenumi}{\textit{(\roman{enumi})}}}
\newcommand{\leialph}{\renewcommand{\labelenumi}{(\alph{enumi})}}
\newcommand{\leiarab}{\renewcommand{\labelenumi}{(\arabic{enumi})}}

\newcommand{\leiialph}{\renewcommand{\labelenumii}{(\alph{enumii})}}
\newcommand{\leiiarab}{\renewcommand{\labelenumii}{(\arabic{enumii})}}

\newenvironment{enumerateroman}
{\leirom \begin{enumerate}}
{\end{enumerate} \leiarab}

\newenvironment{itemizeminus}
  {\begin{itemize}}{\end{itemize}}

\newenvironment{enumeratenumeric}{\begin{enumerate}}{\end{enumerate}}

\newenvironment{enumeratealpha}
{\leialph \leiialph \begin{enumerate}}
{\end{enumerate} \leiarab \leiiarab}

\newenvironment{theproof}{
  \noindent\textbf{Proof.}\ }{\hspace*{\fill}
\begin{math}
\Box
\end{math}
\medskip}

\newenvironment{theproofof}{
  \noindent\textbf{Proof of}\ }{\hspace*{\fill}
\begin{math}
\Box
\end{math}
\medskip}

\newcommand{\tmstrong}[1]{\textbf{#1}}
\newcommand{\tmscript}[1]{\text{\scriptsize $#1$}}

\newtheorem{definition}{Definition}[subsection]
\newtheorem{proposition}{Proposition}[subsection]
\newtheorem{theorem}{Theorem}[subsection]
\newtheorem{varremark}{Remark}[subsection]
\newcommand{\erem}{\hspace*{\fill}$\bigstar$}
\newenvironment{theremark}{\begin{varremark}\em}{\em\erem\end{varremark}}
\newtheorem{lemma}{Lemma}[subsection]

\newcommand{\longrightarrowlim}{\mathop{\longrightarrow}\limits}

\newtheorem{corollary}{Corollary}[subsection]

\newcommand{\nin}{\not\in}

\newcommand{\mathbbm}[1]{\mathbb{#1}}

\newdir{ >}{{}*!/-10pt/\dir{>}}

\newcommand{\A}{\mathbb{A}}
\newcommand{\B}{\mathbb{B}}

\newcommand{\N}{\mathbb{N}}

\newcommand{\M}{\mathbb{M}}
\newcommand{\I}{\mathbb{I}}
\newcommand{\C}{\mathbb{C}}

\newcommand{\R}{\mathbb{R}}

\newcommand{\Sb}{\mathbb{S}}

\newcommand{\ms}{\mathcal{M}}
\newcommand{\as}{\mathcal{A}}
\newcommand{\bs}{\mathcal{B}}
\newcommand{\cs}{\mathcal{C}}
\newcommand{\is}{\mathcal{I}}
\newcommand{\V}{\bs\setminus\as}

\newcommand{\ws}{\mathcal{W}}

\newcommand{\llp}[1]{LLP  \left( #1 \right)}

\newcommand{\rlp}[1]{RLP  \left( #1 \right)}
\newcommand{\free}{\mathbf{\mathcal{F}}}

\newcommand{\forg}{\mathcal{U}}

\newcommand{\deq}{\overset{def.}{=}}
\newcommand{\dee}{\overset{def}{{\Longleftrightarrow}}}

\newcommand{\tcat}{\textup{\textbf{2-Cat}}}
\newcommand{\sset}{\textup{\textbf{sSet}}}
\newcommand{\sett}{\textup{\textbf{Set}}}
\newcommand{\topp}{\textbf{Top}}

\newcommand{\tcatln}{\widetilde{\textup{\textbf{2-Cat}}}}
\newcommand{\cat}{\textup{\textbf{Cat}}}
\newcommand{\tgrph}{\textup{\textbf{2-Grph}}}
\newcommand{\ngrph}{\textup{\textbf{n-Grph}}}

\newcommand{\der}{\textup{\textbf{Der}}}
\newcommand{\lder}{\text{der}}
\newcommand{\id}{\text{id}}
\newcommand{\sd}{\text{Sd}}
\newcommand{\dom}{\text{dom}}
\newcommand{\colim}{\text{colim}}
\newcommand{\cart}{\text{cart}}
\newcommand{\cod}{\text{cod}}
\newcommand{\inn}{\text{in}}
\newcommand{\ex}{\text{Ex}}
\newcommand{\cyl}{\text{Cyl}}
\newcommand{\sesqu}{\textup{\textbf{Sesqu}}}

\newcommand{\nlax}{\textup{\textbf{NLax}}}
\newcommand{\ord}{\textup{\textbf{Ord}}}
\newcommand{\ups}[1]{\uparrow \hspace{-0.25em} #1}
\newcommand{\upa}{{\ups{A}}}

\newcommand{\hc}{{\circ}}
\newcommand{\vc}{{\bullet}}

\newcommand{\presh}[1]{^{#1^{op}}}

\begin{document}

\title{A Model Structure \`a la Thomason on {\tcat}} \author{K. Worytkiewicz, K. Hess, P.E.Parent,A.Tonks} \maketitle

\begin{abstract}
  We exhibit a model structure on {\tcat}, obtained by transfer from
  $\sset$ across the adjunction $C_2 \circ Sd^2 \hspace{0.25em}
  \dashv \hspace{0.25em} Ex^2 \circ N_2$.
\end{abstract}

\section{Introduction}

There are two well-known model category structures on the category $\bold {Cat}$ of small categories: the ``folklore" structure, the existence of which was intuited for many years before it was finally established rigorously by Joyal and Tierney in 1991 \cite{jo-ti}, and the ``topological" structure, developed by Thomason in 1980 \cite {thoma-model} and recently  corrected by Cisinski \cite {cisinski-thoma,cisinski-thesis}.  In the ``folklore" structure, weak equivalences are equivalences of categories, corresponding to a purely category-theoretic view of the role of categories.  On the other hand, the ``topological" structure is defined so that the functor ${Ex^2\circ N:\bold {Cat}\longrightarrow \bold{sSset}}$ induces an equivalence of homotopy categories, where $Ex$ is the right adjoint to the subdivision functor $Sd$, $N$ is the nerve functor and $\bold {sSet}$ is the category of simplicial sets.   In particular, a functor ${F:\Bbb A \longrightarrow \Bbb B}$ between small categories is a weak equivalence if and only if ${NF:N\Bbb A \longrightarrow N\Bbb B}$ is a weak equivalence of simplicial sets.

Our goal in this article is to establish the existence of a Thomason-type, ``topological" model category structure on  $\tcat$, the category of small $2$-categories, complementing Lack's recent proof of the existence of a ``folklore" structure on $\bold {Cat}$ \cite{slack}.    More precisely, we prove that there is a model category structure on $\tcat$ such that  ${Ex^2\circ N_2:\tcat\longrightarrow \bold{sSet}}$ induces an equivalence of homotopy categories, where $N_2$ denotes the $2$-nerve functor. Our methods are analogous to  those of Thomason and Cisinski, though the generalization to $2$-categories is highly nontrivial.

We begin this article with a thorough primer on $2$-category theory in section 2 .  In particular we provide a careful review of the construction of limits and colimits in $\tcat$, as well as of the definition of $N_2$ and its left adjoint, the $2$-categorification functor $C_2$.   We then recall the necessary elements of model category theory in section 3, including a very useful ``Creation Proposition", giving conditions under which model category structure can be transfered across a pair of adjoint functors. 
 
In section 4 we prove the existence of the Thomason-type model category structure on $\tcat$.  We first introduce the notion of right and left ideals of $2$-categories, which we use then in the crucial definitions of a {\sl distortion} between $2$-functors and of a {\sl skew immersion} of $2$-categories.   A distortion from a $2$-functor $F$ to a $2$-functor $G$ is a  sort of left homotopy from $F$ to $G$, which, in fact, induces a simplicial homotopy from $N_2F$ to $N_2G$.  On the other hand, a skew immersion is an inclusion of a left ideal $\as\hookrightarrow \bs$ such that $\as$ is a sort of ``strong deformation retract" (notion defined using distortions) of a right ideal $\ws$ of $\bs$, implying that $N_2\as$ truly is a strong deformation retract  $N_2\ws$ in the usual sense.  The most important example of a skew immersion for our purposes is $(C_2\circ Sd^2)(i_{k,n})$, where ${i_{k,n}:\Lambda ^k[n]\longrightarrow \Delta [n]}$ is a horn inclusion. We establish furthermore that skew immersions are stable under pushout and that the image under $N_2$ of a pushout of a skew immersion along an arbitrary 2-functor is a weak pushout.  Given these reults, we can finally apply the ``Creation Proposition" to to show that $Ex^2\circ N_2$ creates the desired model category structure on $\tcat$.

In the final section of the paper, we show that B\'enabou's ``$2$-category of cylinders" gives a natural path object construction in $\tcat$.  The desire to establish this result motivated the research in this article, as it has an intriguing application in concurrency theory \cite{getco}.

Given a new and interesting model category structure, it is natural to ask what properties the structure satisfies and how well we can characterize fibrations and cofibrations, as well as fibrant and cofibrant objects.  It turns out that the Thomason-type structure on $\tcat$ is both cellular and proper, as we will establish in a future article.  The proof of properness depends on the observation that all cofibrations in $\tcat$ are retracts of skew immersions, since all cofibrations are retracts of elements of $\;\is-cell$, where $\is=\{C_2\circ Sd^2(\partial \Delta [n])\hookrightarrow C_2\circ Sd^2(\Delta [n])\mid n\geq 0\}$, all elements of which are skew immersions.

\section{\label{sec:2-Categories}2-Categories and 2-Nerves}

\subsection{2-Cat}

\subsubsection{2-Graphs}

\begin{definition}
  Let $\mathbb{A}$ be a category. A preglobular object $A$ in $\mathbb{A}$ is
  a $\mathbb{N}$-indexed sequence

  \begin{center}
    $\xymatrix{
\cdots \;\; A_i \ar@<1ex>[r]^{dom_{i-1}} \ar@<-1ex>[r]_{cod_{i-1}} &
A_{i-1 }\;\; \cdots
}$

  \end{center}
  
  {\noindent}of objects and morphisms subject to the identities
  \[ \begin{array}{lcl}
       \dom_i \circ \dom_{i + 1} & = & \dom_i \circ \cod_{i + 1}\\
       \cod_i \circ \dom_{i + 1} & = & \cod_i \circ \cod_{i + 1}
     \end{array} \]
  $A$ is n-truncated if $i < n$. An $n$-graph is a $n$-truncated preglobular
  set.
\end{definition}

\begin{theremark}
  \label{rem:ngrph-toposes}Since an $n$-graph is just a presheaf, {\ngrph} is
  a topos for each $n \in \mathbb{N}$. In particular, {\ngrph} is complete and
  cocomplete.
 \end{theremark}

\begin{definition}
  \begin{enumerateroman}
    \item A graph is a $1$-graph with $\dom \deq \dom_0$ and $\cod \deq
    \cod_0$. Let $H$ be a graph and $a, b \in H_0$, then
    \[ H \left( a, b \right) \deq \left\{ u \in H_1 \hspace{0.25em} \mid
       \hspace{0.25em} \dom \left( u \right) = a \hspace{0.25em} \wedge
       \hspace{0.25em} \cod \left( u \right) = b \right\}  \]
    \item let $G$ be a $2$-graph. As in the case of graphs, the elements of
    $G_0$ are called \textit{vertices} or \textit{$0$-objects} and those of
    $G_1$ \textit{arrows}, \textit{edges} or \textit{$1$-morphisms}. The
    elements of $G_2$ are called \textit{2-cells} or \textit{$2$-morphisms}.
    $G$'s underlying graph $\left\lfloor G \right\rfloor$ is given by its
    1-truncation $G_1 \rightrightarrows G_0$;
    
    \item given $x, y \in G_0$, $G \left( x, y \right)$is the graph with
    \[ \begin{array}{lcl}
         G \left( x, y \right)_0 & \deq & \left\{ f \in G_1 \hspace{0.25em}
         \mid \hspace{0.25em} \dom_0 \left( f \right) = x \hspace{0.25em}
         \wedge \hspace{0.25em} \cod_0 \left( f \right) = y \right\}\\
         G \left( x, y \right)_1 & \deq & \left\{ \alpha \in G_2
         \hspace{0.25em} \mid \hspace{0.25em} \dom_1 \left( \alpha \right),
         \hspace{0.25em} \cod_1 \left( \alpha \right) \in G \left( x, y
         \right)_0 \right\}
       \end{array} \]
    and with $\dom_{x, y} \hspace{0.25em}, \cod_{x, y} : \hspace{0.25em} G
    \left( x, y \right)_1 \rightarrow G \left( x, y \right)_1$ given by
    \[ \begin{array}{lcl}
         \dom_{x, y} \left( \alpha \right) & \deq & \dom_1 \left( \alpha
         \right)\\
         \cod_{x, y} \left( \alpha \right) & \deq & \cod_1 \left( \alpha
         \right)
       \end{array} \]
  \end{enumerateroman}
\end{definition}

Properties and concepts defined with respect to $G \left( x, y
\right)$ (or its more structured counterparts to be introduced below) are
called \textit{local}. For instance, a morphism of graphs $h : \hspace{0.25em}
G \rightarrow H$ is \textit{locally injective if $h_1 \mid_{G \left( x, y
\right)}$} is an injective function for each $x, y \in G_0$.

\subsubsection{Derivation Schemes and Sesquicategories}

\begin{definition}
  A derivation scheme is a 2-graph $D$ such that the underlying graph
  $\left\lfloor D \right\rfloor$ is a category. The composition in
  $\left\lfloor D \right\rfloor$ is denoted $\hc$ and written infix in the
  evaluation order. Morphisms of derivation schemes are morphisms of 2-graphs
  that are functors on the underlying categories.
\end{definition}

\begin{proposition}
  Derivation schemes and their morphisms form the category {\der}. There is an
  adjunction
  
  \begin{center}
    $\xymatrix{ 
\der \rtwocell^{F_{\mathbf \der}}_{U_{\mathbf \der}}{`{\perp}} & \tgrph
}$

  \end{center}

\end{proposition}
\begin{theproof}
Let $G$ be a $2$-graph. The free derivation scheme $F_{\lder} ( G )$ is given
by
\[ \left\lfloor F_{\lder} \left( G \right) \right\rfloor = \mathcal{F} \left(
   \left\lfloor G \right\rfloor \right) \]
where $\mathcal{F} \left(\left\lfloor G \right\rfloor \right)$ is the free category on $\left\lfloor G \right\rfloor$.
\end{theproof}

Let $x, y \in G_0$. A situation involving an $\alpha \in G \left( x, y
\right)_1$ such that $\dom \left( \alpha \right) = f$ and $\cod \left( \alpha
\right) = g$ is customarily drawn as

\begin{center}
  $\xymatrix{ x \rtwocell^f_g{\alpha} & y 
}$

\end{center}

\begin{definition}
  \label{def:A-sesquicategory}A sesquicategory $\mathbb{S}$ is a derivation
  scheme such that $\mathbb{S} \left( x, y \right)$ is a category for all $x,
  y \in \mathbb{S}_0$. The composition in $\mathbb{S} \left( x, y \right)$ is
  denoted $\vc$ and is written infix in the evaluation order. For each $x', x,
  y \in \mathbb{S}_0$ there is an operation
  \[ W_{\tmop{left}} : \hspace{0.25em} \mathbb{S} \left( x', x \right)_0
     \times \mathbb{S} \left( x, y \right)_1 \rightarrow \mathbb{S} \left( x',
     y \right)_1 \]
  and for each $x, y, y' \in \mathbb{S}_0$ there is an operation
  \[ W_{\tmop{right}} : \hspace{0.25em} \mathbb{S} \left( x, y \right)_1
     \times \mathbb{S} \left( y, y' \right)_0 \rightarrow \mathbb{S} \left( x,
     y' \right)_1 \]
  Both operat{\noindent}ions are called \textit{whiskering} and are denoted
  $\hc$ by abuse of notation. $W_{\tmop{left}}$ is subject to the identities
  \begin{enumerate}
    \item given
    
    \begin{center}
      $\xymatrix{ x \ar[r]^{\id} & x \rtwocell^f_g{\alpha} & y 
}$

    \end{center}
    
    {\noindent}the equation
    \[ \alpha \hc \id_x = \alpha \]
    holds;
    
    \item given
    
    \begin{center}
      $\xymatrix{ x' \ar[r]^{f} & x \rtwocell^u_u{\id} & y
}$

    \end{center}
    
    {\noindent}the equation
    \[ \id_u \hc f = \id_{u \circ f} \]
    holds;
    
    \item given
    
    \begin{center}
      $\xymatrix{ x'' \ar[r]^{f'} & x' \ar[r]^{f} & x \rtwocell^u_u{\alpha} 
& y
}$

    \end{center}
    
    {\noindent}the equation
    \[ \alpha \hc \left( f \hc f' \right) = \left( \alpha \hc f \right) \hc f'
    \]
    holds;
    
    \item given
    
    \begin{center}
      $\xymatrix{ x' \ar[r]^{f} & x \ruppertwocell^{u}{\alpha} \rlowertwocell_{w}{\beta} \rto_(.35)v & y 
}$

    \end{center}
    
    {\noindent}the equation
    \[ \left( \beta \vc \alpha \right) \hc f = \left( \beta \hc f \right) \vc
       \left( \alpha \hc f \right) \]
    holds;
    
    \item the rules governing $W_{right}$ are defined symmetrically;
    
    \item given
    
    \begin{center}
      $\xymatrix{ x \ar[r]^{f} & x \rtwocell^f_g{\alpha} & y \ar[r]^{g} & z
}$

    \end{center}
    
    {\noindent}the equation $g \circ \left( \alpha \circ f \right) = \left( g
    \circ \alpha \right) \circ f$ holds.
  \end{enumerate}
  Morphisms of sesquicategories, called \textit{sesqifunctors,} are
  \textit{}morphisms of the underlying derivation schemes which are locally
  functors and which preserve whiskering.
\end{definition}

The equations of a sesquicategory guarantee in particular that there is no
harm to write the 2-cells as strings like
\[ g_m \circ \cdots \circ g_1 \circ \alpha \circ f_n \cdots f_1 \]
\begin{proposition}
  Sesquicategories and sesquifunctors form in the category {\sesqu} .
  There is an adjunction
  \begin{center}
    \begin{center}
      $\xymatrix{
\sesqu \rtwocell^{F_{\mathbf \sesqu}}_{U_{\mathbf \sesqu}}{`{\perp}} & \der
}$

    \end{center}
  \end{center}
\end{proposition}

A free sesquicategory $\mathcal{F} \mathbb{D}$ over a derivation scheme
$\mathbb{D}$ is given by formally adding all the whiskering composites and all
the vertical composites.

\begin{definition}
  Let $\Sb$ be a sesquicategory. A \textit{sesquicongruence} on $\mathbb{S}$
  is a family
  \[ \left\{ \sim^1_{X, Y} \subseteq \mathbb{A} \left( X, Y \right) \times
     \mathbb{A} \left( X, Y \right) \right\}_{X, Y \in \mathbb{A}_0} \]
  of equivalence relations on morphisms and a family
  \[ \left\{ \sim^2_{f, g} \subseteq \mathbb{A} \left( X, Y \right) ( f, g )
     \times \mathbb{A} \left( X, Y \right) ( f, g )
     \right\}_{\tmscript{\begin{array}{l}
       X, Y \in \mathbb{\Sb}_0\\
       f, g \in \Sb ( X, Y )
     \end{array} }} \]
  of equivalence relations on 2-cells such that
  \begin{enumerateroman}
    \item $\alpha \sim^2 \beta \Longrightarrow \theta \bullet \alpha \bullet
    \varphi \sim^2 \theta \bullet \alpha \bullet \varphi \text{ and } g \circ
    \alpha \circ f \sim^2 g \circ \beta \circ f$
    
    \item $f \sim^1 g \Longrightarrow \phi \circ f \circ \psi \sim^2 \phi
    \circ g \circ \psi$ 
    
    \item $\id_f \sim^2  \id_g \Longrightarrow f \sim^1 g$
  \end{enumerateroman}
\end{definition}

\begin{theremark}
  In particular, $\sim^1$ is a congruence on $\left\lfloor \Sb \right\rfloor$.
\end{theremark}

\begin{proposition}
  An arbitrary intersection of sesquicongruences is again a sesquicongruence.
  The quotient $\mathbb{} \Sb / \sim$ of a sesquicategory $\Sb$ by a
  sesquicongruence $\sim$ is again a sesquicategory.
\end{proposition}

\subsubsection{2-Categories}

\begin{definition}
  Let $\mathbb{S}$ be a sesquicategory and $x, y, z \in \mathbb{S}_0$. The
  latter satisfy the \textit{interchange law} if any diagram of the form
  
  \begin{center}
    $\xymatrix{ x \rtwocell^f_g{\alpha} & y \rtwocell^{f'}_{g'}{\alpha'} & z
}$

  \end{center}
  
  {\noindent}verifies the equation
  \[ \left( g^{\prime} \hc \alpha \right) \vc \left( \alpha^{\prime} \hc f
     \right) = \left( f^{\prime} \hc \alpha \right) \vc \left( \alpha^{\prime}
     \hc g \right) \hspace{2em} \left( \ast \right) \]
  A 2-category is a sesquicategory in which the interchange law holds for
  every triple of objects. A 2-functor is a sesquifunctor between 2-categories.
  2-categories and 2-functors form the category  {\tcat}. 
\end{definition}

\begin{theremark}
  The quotient of a 2-category by a sesquicongruence is again a 2-category.
\end{theremark}

\begin{proposition}
  \textbf{\label{pro:Gray}} {\dueto{Gray {\cite{gray:2-cat}}}} The functor
  $\mathcal{\left\lfloor_- \right\rfloor} : \tcat \longrightarrow \cat$ which forgets the 2-cells
  has a right adjoint.
\end{proposition}

\begin{theproof}
  The right adjoint turns a homset into a trivial connected groupoid.
\end{theproof}

The interchange law is often called by the name of R.Godement {\cite{god}}.
A 2-category $\mathcal{A}$ admits in particular a ``horizontal'' composition
of 2-cells where $\alpha' \hc \alpha$ is given by either side of $\left( \ast
\right)$, giving rise to a family of functors
\[ \_ \hc \_ : \hspace{0.25em} \mathcal{A} \left( y, z \right) \times
   \mathcal{A} \left( x, y \right) \rightarrow \mathcal{A} \left( x, z \right)
\]
indexed by triples $x, y, z \in \mathcal{A}_0$. This is the way 2-categories
are usually introduced in the literature (c.f. {\cite{bor1}}), while the
exposition above is drawn from {\cite{street-cat-struct}}.

\begin{proposition}
  There is an adjunction
  
  \begin{center}
    $\xymatrix{ 
\tcat \rtwocell^{F_{\mathbf \tcat}}_{U_{\mathbf \tcat}}{`{\perp}} & \sesqu
}$

  \end{center}
\end{proposition}

It is easy to see that constructing the free 2-category on a sesquicategory
amounts to quotienting the latter by the sesquicongruence generated by the
equations enforcing the Godement law for all triples of objects. We thus have
the series of adjunctions

\begin{center}
  $\xymatrix{ 
\tcat \rtwocell^{F_{\mathbf \tcat}}_{U_{\mathbf \tcat}}{`{\perp}} 
& 
\sesqu \rtwocell^{F_{\mathbf \sesqu}}_{U_{\mathbf \sesqu}}{`{\perp}} 
&
\der \rtwocell^{F_{\mathbf \der}}_{U_{\mathbf \der}}{`{\perp}} 
& 
\tgrph
}$

\end{center}

\begin{definition}
  Let $G$ be a 2-graph and
  \[ \free  \deq F_{\tmmathbf{\tcat}} \circ F_{\tmmathbf{\sesqu}} \circ
     F_{\tmmathbf{\der}} \]
  The \textit{free 2-category $\mathcal{F} ( G )$ on $G$} is given by this
  functor.
\end{definition}

A free 2-category on a 2-graph (or a derivation scheme) involves thus
``horizontal'' sequences in dimension 1 and 2 as well as ``vertical''
sequences in dimension 2. We write
\[ < f_1 ; \cdots ; f_n > \]
for a 1-dimensional horizontal sequence of morphisms,
\[ \ll A_1 ; \cdots ; A_n \gg \]
for a horizontal sequence of morphisms and/or 2-cells and
\[ \ll \alpha_1 : \cdots : \alpha_m \gg \]
for a vertical sequence of 2-cells. We define the concatenation operations
\[ \text{$\ll A_1 ; \cdots ; A_k \gg ; \ll A_{k + 1} ; \cdots ; A_n \gg = \ll
   A_1 ; \cdots ; A_n \gg$} \]
and
\[ \ll \alpha_1 : \cdots : \alpha_l \gg : \ll \alpha_{l + 1} : \cdots :
   \alpha_m \gg = \ll \alpha_1 : \cdots : \alpha_m \gg \]
at any index. Those are obviously associative and can be mixed whenever it
makes sense, e.g.
\[ \ll \alpha : \alpha' \gg ; \ll \beta : \beta' \gg = \ll \alpha ; \beta \gg
   : \ll \alpha' ; \beta' \gg \]
is an instance of the interchange law. Domains and codomains are usually clear
from context. If not, we indicate them as subscripts. In case of
endomorphisms or endo-2-cells we do not duplicate those subscripts, e.g
\[ < >_X \]
is the empty sequence with domain and codomain $X$, i.e. the 1-dimensional
identity at $X$. Similarly,
\[ \ll \gg_f \]
is the 2-dimensional identity at $f$.

\subsection{Limits and Colimits in {\tcat}}

\begin{proposition}
  \label{pro:2cat-complete} {\tcat} is complete and cocomplete.
\end{proposition}

\begin{theproof}
  Limits are obvious. Let $D : \I \longrightarrow \text{{\tcat}}$ be a
  diagram. There is the colimiting cocone
  \[ \left\{ \iota_K : ( \forg \circ D ) \longrightarrow \colim ( \forg \circ
     D ) \right\}_{K \in \I} \]
  in {\tgrph}. Consider the 2-category
  \[ \free \left( \colim ( \forg \circ D ) \right) / \sim \]
  where $\sim$ is the sesquicongruence generated by
  \begin{enumerateroman}
    \item $\ll \iota_K ( \alpha ) : \iota_K ( \beta ) \gg = \ll \iota_K (
    \beta \bullet \alpha ) \gg$
    
    \item $\ll \iota_K ( \id_f ) \gg = \ll \gg_f$
    
    \item $\ll \iota_K ( f ) ; \iota_K ( \alpha ) ; \iota_K ( g ) \gg = \ll
    \iota_K ( g \circ \alpha \circ f ) \gg$
    
    \item $< \iota_K ( \id_X ) > = < >_X$ 
  \end{enumerateroman}
  for all $K \in \I$. There is the cocone
  \[ \left\{ \kappa_K : D ( K ) \longrightarrow \free \left( \colim ( \forg
     \circ D ) \right) / \sim \right\}_{K \in \I} \]
  in {\tcat}, given by
  \[ \left( \kappa_K \right)_2 ( \alpha ) \deq \ll \iota_K ( \alpha ) \gg \]
  and
  \[ \left( \kappa_K \right)_1 ( f ) \deq < \iota_K ( f ) > \]
  This cocone is colimiting. To see this, suppose there is a cocone
  \[ \left\{ c_K : D ( K ) \longrightarrow \cs \right\}_{K \in \I} \]
  over $D$. Then there is the comparison morphism
  \[ m : \colim ( \forg \circ D ) \longrightarrow \forg ( \cs ) \]
  in {\tgrph}. Its transpose
  \[ \bar{m} : \free \left( \colim ( \forg \circ D ) \right) \longrightarrow
     \cs \]
  over the adjunction $\free \dashv \forg$ remains defined after the passage
  to the quotient and is the desired comparison morphism.
\end{theproof}

Our proof above, one of the manifold possible variants, generalizes Gabriel's
and Zisman's construction of colimits in $\cat$ (c.f. {\cite{GZ67}}). It is
easy to see that our construction amounts to doing first the construction 
on the
underlying category as in {\cite{GZ67}} and then to taking care of the 2-cells.
It {\tmem{has}} to be that way because of proposition \ref{pro:Gray}.

\begin{theremark}
  \label{rem-pushinc}The calculatory recipe given in the proof of proposition
  \ref{pro:2cat-complete} is quite practical indeed. Consider for instance the
  case of pushing inclusions out:

  \begin{center}
    $
\xy
\xymatrix{
\as \ar[r]^F \ar@{ >->}[d] & \as' \ar@{ >->}[d]^{\kappa_{\as'}}
\\
\ws \ar[r]_(.4){\kappa_{\ws}} & \as' +_\as \ws
}
\POS( 14 , -8);\POS( 14,-10.4 ) \connect@{-}
\POS( 14 , -8);\POS( 16.4,-8 ) \connect@{-}
\endxy
$

  \end{center}

  {\noindent}Then there is the pushout square

  \begin{center}
    $
\xy
\xymatrix{
\forg (\as) \ar[r]^{\forg (F)} \ar@{ >->}[d]
& 
\forg (\as') \ar[d]^{\iota_{\as'}}
\\
\forg (\ws) \ar[r]_{\iota_{\ws}} & P
}
\POS( 14.5 , -9);\POS( 14.5,-11.4 ) \connect@{-}
\POS( 14.5 , -9);\POS( 16.9,-9 ) \connect@{-}
\endxy
$

  \end{center}
  
  {\noindent}in {\tgrph}, where
  \[ P \cong \left( \as'_i + \left( \ws_i \backslash \as_i \right) \right)_{0
     \leqslant i \leqslant 2} \]
  with structural maps given by universal property as the copairs
  \[ \partial^i_P = \left[ \inn_{\as_i'} \circ \partial^i_{\as'}, \left( F_i +
     \id_{\ws_i \backslash \as_i} \right) \circ \partial^i_{\ws} |_{_{\ws_{i +
     1} \backslash \as_{i + 1}}} \right] \]
  for $i \in \{ 0, 1 \}$ and $\partial \in \{ \dom, \cod \}$. On the other
  hand
  \[ \iota_{\as'} = \left( \inn_{\as'_2}, \inn_{\as'_1}, \inn_{\as'_0} \right)
  \]
  and
  \[ \iota_{\ws} = \left( F_2 + \id_{\ws_2 \backslash \as_2}, F_1 + \id_{\ws_1
     \backslash \as_1}, F_0 + \id_{\ws_0 \backslash \as_0} \right) \]
  Then
  \begin{eqnarray*}
    \left( \kappa_{\as'} \right)_0 & = & \inn_{\as'_0}\\
    \left( \kappa_{\as'} \right)_1 ( f ) & = & < \inn_{\as'_1} ( f ) >\\
    \left( \kappa_{\as'} \right)_2 ( \alpha ) & = & \ll \inn_{\as'_2} ( \alpha
    ) \gg
  \end{eqnarray*}
  determines a morphism of 2-graphs $\kappa_{\as' } : \as' \longrightarrow
  \free ( P )$ while
  \begin{eqnarray*}
    \left( \kappa_{\ws} \right)_0 & = & F_0 + \id_{\ws_0 \backslash \as_0}\\
    \left( \kappa_{\ws} \right)_1 ( u ) & = & < \left( F_1 + \id_{\ws_1
    \backslash \as_1} \right) ( u ) >\\
    \left( \kappa_{\ws} \right)_2 ( \theta ) & = & < \left( F_2 + \id_{\ws_2
    \backslash \as_2} \right) ( \theta ) >
  \end{eqnarray*}
determines a morphism of 2-graphs
  \[ \kappa_{\ws} : \ws \longrightarrow \free ( P )\] so
  \[ \as' +_{\as}  \ws \cong \free ( P ) / \sim \]
  with $\sim$ the smallest sesquicongruence making $\kappa_{\as'}$ and
  $\kappa_{\ws}$ 2-functorial.
  
  It follows that inclusions in {\tcat} are stable under pushout. In
  particular, if an inclusion is full and locally full, then pushing it out
  will result in a full and locally full one. 
\end{theremark}

\subsection{2-Nerve and 2-Categorification}

\subsubsection{Simplicial Sets}

\begin{lemma}
  {\dueto{Kan}} \label{lem:kan}Let $F : \hspace{0.25em} \mathbb{C} \rightarrow
  \mathbb{A}$ be a functor and $A \in \mathbb{A}$. The assignment
  \[ A \mapsto \mathbb{A} \left( F \left( \_ \right), A \right) \]
  determines a functor $F_{\ast} : \hspace{0.25em} \mathbb{A} \rightarrow
  \sett^{\mathbb{C}^{op}}$. If $\mathbb{A}$ is cocomplete then $F_{\ast}$ has
  a left adjoint $F_!= \tmop{Lan}_y F$ and $F$ factors through $F_!$ by the Yoneda embedding
  $y : \hspace{0.25em} \mathbb{C} \rightarrow \sett^{\mathbb{C}^{op}}$:
\begin{center}
    $\xymatrix{ 
\presh{\C} \drtwocell<2>^{F_!}_{F_*}{'{\vdash}} & \\
\C \ar[u]^y \ar[r]_F & \A 
}$

  \end{center}
\end{lemma}

{\noindent}The condition of $\mathbb{A}$ being cocomplete is
stronger than the existence of the relevant Kan extension, yet 
it is verified in most of the cases of interest.

\begin{definition}
  Let $[ n ] \deq \{ 0 < 1 \cdots < n \}$ be the $n$th finite ordinal and
  \begin{itemizeminus}
    \item $\delta^i_n : [ n - 1 ] \longrightarrow [ n ]$ be the increasing
    injection missing $i$;
    
    \item $\sigma^i_n : [ n + 1 ] \longrightarrow [ n ]$ be the non decreasing
    surjection taking twice the value $i$;
  \end{itemizeminus}
  The category $\Delta$ has finite ordinals as objects and is generated by
  \[ \{ \delta^i_n | n \in \N, 0 < n, 0 \leqslant i \leqslant n \} \cup \{
     \sigma^i_n | n \in \N, 0 \leqslant i \leqslant n \} \]
  Let $\C$ be a category. The category of {\tmem{simplicial objects}} in
  $\C$ is $\C^{\Delta^{\tmop{op}}}$ while the category of {\tmem{cosimplicial
  objects}} in $\C$ is $\C^{\Delta}$.  
\end{definition}

As a matter of terminology, if the objects of $\C$ are called ``gadgets'' then
(co)simplicial objects in $\C$ are called ``(co)simplicial gadgets'', e.g.
simplicial sets, simplicial groups, simplicial 2-categories and so on. It is
customary to write {\sset} for the category of simplicial sets and, given $K
\in \tmmathbf{\sset}$, to abbreviate $K_n  \deq K ( [ n ] )$.

\begin{definition}
  Let $K \in \tmmathbf{\sset}$. An element of $K_n$ is called an
  {\tmem{$n$-simplex}}. The representable prefsheaf $\Delta [ n ] \deq \Delta
  (_-, [ n ] ) \in \tmmathbf{\sset}$ is called the {\tmem{standard
  $n$-simplex}}. An $n$-simplex is a {\tmem{face}} if it is in the image of
  some $\partial_i^n  \deq K ( \delta_i^n )$. It is {\tmem{degenerate}} if it
  is in the image of some $\varepsilon_i^n  \deq K ( \sigma_n^i )$. A
  simplicial set is $n$-skeletal if the $m$-simplices are degenerate for $m >
  n$.
\end{definition}

\begin{theremark}
  The standard $n$-simplex $\Delta [ n ]$ is $n$-skeletal. It has precisely
  one non-degenerate $n$-simplex, namely $\tmop{id}_{[ n ]} \in \Delta ( [ n
  ], [ n ] )$. The other degenerate $m$-simplices are all faces.
\end{theremark}

\begin{definition}
  \begin{enumerateroman}
    \item The subobject $\partial \Delta [ n ] \rightarrowtail \Delta [ n ]$,
    obtained from $\Delta [ n ]$ by removing $\tmop{id}_{[ n ]}$, is called
    boundary;
    
    \item Let $1 \leqslant k \leqslant n + 1$; the subobject $\Lambda^k [ n ]
    \rightarrowtail \Delta [ n ]$, obtained from $\partial \Delta [ n ]$ by
    removing $\partial_k^n ( \tmop{id}_{[ n ]} ) = \delta^k_n$, is called
    $k$th {\tmem{horn}}.
  \end{enumerateroman}
\end{definition}

\begin{proposition}
  \label{def-georel}Let
  \[ \Delta_n  \deq  \left\{ ( t_1, \ldots, t_n ) \in \R^n | \sum_{i = 1}^n
     t_i = 1 \wedge \forall 1 \leqslant i \leqslant n. t_i \geqslant 0
     \right\} \]
  be the standard topological $n$-simplex. The functor
  \[ \begin{array}{llll}
       g : & \Delta & \longrightarrow & \tmmathbf{\topp}\\
       & [ n ] & \longmapsto & \Delta_n
     \end{array} \]
  determines an adjunction
  \[ g_! = |_- | \dashv \tmop{Sing} = g_*\]
\end{proposition}

The left adjoint gives the geometric realization of a simplicial set while the
right adoint gives the singular complex of a topological space. In particular,
singular homology is a special case of simplicial homology via this right
adjoint.

\subsubsection{Orientals}

\begin{definition}
  Let $\left[ n \right] \in \Delta$, $\delta_{i, j}$ be the inequality $i
  \leqslant j$ seen as a morphism in $[ n ]$ and $\bar{\Delta}_n$ be the
  derivation scheme given by the data
  \begin{enumerate}
    \item $\left| \bar{\Delta}_n \right| \deq \mathcal{F} \left( \left[ n
    \right] \right)$;
    
    \item $\left(\bar{\Delta}_n\right)_2 \deq \left\{ \delta_{i, j, k} \mid
    \hspace{0.25em} 0 \leq i < j < k \leq n \right\}$ where
    \[ \dom_1 \left( \delta_{i, j, k} \right) = < \delta_{i, j} ; \delta_{j,
       k} > \]
    and
    \[ \cod_1 \left( \delta_{i, j, k} \right) = < \delta_{i, k} > \]
  \end{enumerate}
  The 2-category $\Delta_n$ is the free 2-category $\mathcal{F} \left(
  \bar{\Delta}_n \right)$ over $\bar{\Delta}_n$ quotiented by the relations
  \[ < \delta_{i, j, k} ; \delta_{k, l} > : \ll \delta_{i, k, l}
     \hspace{0.25em} \gg = \hspace{0.25em} < \delta_{i, j} ; \delta_{j, k, l}
     > : \ll \delta_{j, k, l} \gg \]
\end{definition}

Following  Street {\cite{street-osimpl}}, we call the $\Delta_n$'s
\textit{2-orientals}.

\begin{proposition}
  The construction $\Delta_{(_- )} : \hspace{0.25em} \Delta \longrightarrow
  \tmmathbf{\tcat}$ is functorial and determines an adjunction
  \[ C_2 \dashv N_2 \]
\end{proposition}

\begin{theproof}
  The functoriality is immediate while $C_2 \deq {\Delta_{(\_)}}_!$ and $N_2
  \deq {\Delta_{(\_)}}_{\ast}$.
\end{theproof}

We call $N_2$ \textit{2-nerve} and $C_2$ \textit{2-categorification}.

\begin{theremark}
  Given a simplicial set $K$, $C_2 \left( K \right)$ is the free 2-category on
  the derivation scheme determined by $\left( K_i \right)_{0 \leq i \leq 2}$,
  quotiented by the sesquicongruence generated by $K_3$.
\end{theremark}

\subsection{Normal Lax Functors}

\begin{definition}
  Let $\mathcal{A}$ and $\mathcal{B}$ be 2-categories and $F : \hspace{0.25em}
  \mathcal{A} \rightarrow \mathcal{B}$ a morphism of the underlying 2-graphs.
  $F$ is a \textit{normal lax functor} provided
  \begin{enumerateroman}
    \item it is locally a functor;
    
    \item it preserves horizontal identites;
    
    \item for any $f \in \mathcal{A} \left( x, y \right)$ and $g \in
    \mathcal{A} \left( y, z \right)$ there is the structural 2-cell
    \[ \gamma_{f, g} : \hspace{0.25em} F \left( g \right) \hc F \left( f
       \right) \Rightarrow F \left( g \hc f \right) \]
    such that
    \begin{enumeratealpha}
      \item given any $h \in \mathcal{A} \left( z, a \right)$, the equation
      \[ \gamma_{g \hc f,h} \vc \left( F \left( h \right) \hc \gamma_{f, g}
         \right) = \gamma_{f, h \hc g} \vc \left( \gamma_{g, h} \hc F \left( f
         \right) \right) \]
      holds;
      
      \item given any $\alpha : f \Rightarrow f'$ and $\beta : g \Rightarrow
      g'$, the equation
      \[ \gamma_{f', g'} \bullet \left( F ( \beta ) \circ F ( \alpha ) \right)
         = F ( \beta \circ \alpha ) \bullet \gamma_{f, g} \]
      holds.
    \end{enumeratealpha}
  \end{enumerateroman}
\end{definition}

\begin{theremark}
  A 2-functor is thus a special case of a normal lax functor
 where the structural
  2-cells are all identities.
\end{theremark}

\begin{theremark}
  \label{rem-2catlax}Normal lax functors compose in the obvious way and 
  this composition is associative. The category of 2-categories and normal lax
  functors $\tcatln$ has the usual products, yet it is not finitely complete. 
\end{theremark}

\begin{theremark}
  \label{rem:2-nerve}Let $\tmmathbf{\nlax} \left( \left[ n \right],
  \mathcal{A} \right)$be the set of normal lax functors from $\left[ n
  \right]$ to $\mathcal{A}$. Then
  \[ N_2 \left( \mathcal{A} \right)_n = \tmmathbf{\nlax} \left( \left[ n
     \right], \mathcal{A} \right) \]
  and $N_2$ acts on 2-functors by postcomposition.
  Let $K$ be a simplicial set and let us write $S_{_{i_0, \ldots, i_n}} \in
  K_n$ where $i_0 < \cdots < i_n$ for an $n$-simplex. We use the notation
  \[ \partial_j ( S_{i_0, \ldots, i_n} ) \deq S_{i_0, \ldots, i_{j - 1}, i_{j
     + 1}, \ldots, i_n} \]
  for the faces. The assignment
  \[ \delta_{p, q} \longmapsto < S_{i_p, i_q} > \]
  determines a normal lax functor $S : [ n ] \longrightarrow C_2 ( K )$ with
  the structural 2-cells
  \[ \gamma_{p, q, r} = \ll S_{i_p, i_q, i_r} \gg \]
  The unit $\eta_K : K \longrightarrow ( N_2 \circ C_2 ) ( K )$ of the
  adjunction $C_2 \dashv N_2$ is the simplicial map given in degree $n$ by
  \[ S_{_{i_0, \ldots, i_n}} \longmapsto S \]
\end{theremark}

\begin{theremark}
  \label{rem-n2c2n1}Let $\mathbb{A}$ be a category and $N_1 : \cat \rightarrow
  \sset$ be the usual categorical nerve. Let us write $[ f_1, \ldots, f_n ]$
  for a composable sequence of arrows seen as an $n$-simplex in the nerve.
  $\left( C_2 \circ N_1 \right) \left( \mathbb{A} \right)$ can be
  characterized as follows: the objects are those of $\mathbb{A}$, the arrows
  are generated by those of $\mathbb{A}$ (they are formal composites), while
  the 2-cells are generated by the collection
  \[ [ f, g ] : < f ; g > \Longrightarrow < g \circ f > \]
  subject to the relations
  \[ \ll f ; [ g, h ] \gg : \ll [ f, h \circ g ] \gg = \ll [ f, g ] ; h \gg :
     \ll [ g \circ f, h ] \gg \]
  In particular, $\eta_{N_1 ( \A )}$ is an iso of simplicial sets for any
  category $\A$ by remark \ref{rem:2-nerve}.
\end{theremark}

\section{\label{sec:Model-Category-Theory}Model Category Theory}

In this section, we review some classical and less classical material 
about model categories. Most of the section on topoi is included 
because of its intrinsic beauty.

\subsection{Basic Facts about Model Categories}

\begin{definition}
  Let $\mathbb{M}$ be a category. $\mathcal{L}, \mathcal{R} \subseteq
  \mathbb{M}_1$ form a \textit{weak factorization system $\left( \mathcal{L},
  \mathcal{R} \right)$} if
  \begin{enumerate}
    \item any morphism $f \in \mathbb{M}_1$ factors as $f = r \hc l$ with $r
    \in \mathcal{R}$ and $l \in \mathcal{L}$;
    
    \item $\mathcal{R} = \rlp{\mathcal{L}}$ and $\mathcal{L} =
    \llp{\mathcal{R}}$.
  \end{enumerate}
\end{definition}

\begin{definition}
  $\mathbb{M}$ is a \textit{model category} if it is complete, cocomplete and
  has three distinguished classes of morphisms $\mathcal{C}, \mathcal{W},
  \mathcal{F} \subseteq \mathbb{M}_1$ such that
  \begin{enumerate}
    \item $\left( \mathcal{C}, \mathcal{F} \cap \mathcal{W} \right)$and
    $\left( \mathcal{C} \cap \mathcal{W}, \mathcal{F} \right)$ are weak
    factorization systems;
    
    \item $\mathcal{C}$, $\mathcal{F}$ and $\mathcal{W}$ are closed under
    retracts in $\mathbb{M}^{\rightarrow}$;
    
    \item if two of the morphisms in a commuting triangle are in $\mathcal{W}$
    so is the third one.
  \end{enumerate}
\end{definition}

It is established terminology to call morphisms in $\mathcal{F}$
\textit{fibrations} with $\twoheadrightarrow$ as notation, those in
$\mathcal{C}$ \textit{cofibrations} with $\rightarrowtail$ as notation and
those in $\mathcal{W}$ \textit{weak equivalences} with
$\longrightarrowlim^{\sim}$as notation. It is also customary to call morphisms
in $\mathcal{F} \cap \mathcal{W}$ \textit{acyclic fibrations} and those in
$\mathcal{C} \cap \mathcal{W}$ \textit{acyclic cofibrations}.

\begin{definition}
  Let $\mathbb{M}$ be a cocomplete category and $I \subseteq \mathbb{M}_1$.
  \begin{enumerate}
    \item Let $\lambda$ be an ordinal. A \textit{$\left( \lambda, I
    \right)$-sequence in $\mathbb{M}$} is a cocontinous functor $\lambda
    \rightarrow \mathbb{M}$ such that all its values on morphisms are in $I$.
    
    \item $A \in \mathbb{M}$ is \textit{small with respect to} $I$ if there is
    a cardinal $\kappa$ such that the covariant hom-functor $\mathbb{M} \left(
    A, \_ \right)$ preserves colimits of all $\left( \lambda, I
    \right)$-sequences for all regular cardinals $\lambda \geq \kappa$ .
    
    \item $I$ \textit{permits the small object argument} if the domains of
    morphisms in $I$ are small with respect to $I$.
  \end{enumerate}
\end{definition}

\begin{definition}
  A model category $\mathbb{M}$ is \textit{cofibrantly generated} if there are
  sets of morphisms $I, J \subseteq \mathbb{M}_1$ permitting the small object
  argument and such that
  \[ \mathcal{F} \cap \mathcal{W} = \rlp{I} \textrm{} \]
  and
  \[ \mathcal{F} = \rlp{J} \]
\end{definition}

$I$ is called the set of the generating cofibrations while $J$ is called the
set of generating acyclic cofibrations, this since

\begin{proposition}
  Morphisms in $I$ are cofibrations while those in $J$ are acyclic
  cofibrations.
\end{proposition}

\begin{definition}
  A continous map $f : X \longrightarrow Y$ is a {\tmem{weak homotopy
  equivalence}} if
  \[ \pi_n ( f, x ) : \pi_n ( X, x ) \longrightarrow \pi_n ( Y, f ( x ) ) \]
  is a homeomorphism for any choice of the basepoint $x \in X$.
\end{definition}

\begin{theorem}
  \label{theo-quillen} {\dueto{Quillen}}There is a cofibrantly generated model
  structure on $\tmmathbf{\tmop{Top}}$ such that
  \begin{itemizeminus}
    \item the weak equivalences are the weak homotopy equivalences;
    
    \item $I = \{ S^{n - 1} \hookrightarrow D^n | n \geqslant 0 \}$;
    
    \item $J = \{ I^{n - 1} \times \{ 0 \} \hookrightarrow I^n | n \geqslant 0
    \}$.
  \end{itemizeminus}
\end{theorem}

The model structure of theorem \ref{theo-quillen} is called the ``standard''
or Serre model structure on \textbf{Top}.

\subsection{Model Structures on Topoi of Presheaves}

One of those topoi, namely 
$\sset$, is ubiquitous in homotopy theory:

\begin{theorem}\label{theo-quillen-sset}
  {\dueto{Quillen}} $\sset$ is a cofibrantly generated model category with
  \begin{itemizeminus}
    \item $\ws = \{ f \in \tmmathbf{\sset}_1 | |f| \in
    \ws_{\tmmathbf{\tmop{Top}}} \}$;
    
    \item $\cs = \{ \tmop{Monos} \}$;
    
    \item $I = \left\{ \partial \left[ n \right] \rightarrowtail \Delta \left[
    n \right] | n \in \mathbb{N} \right\}$;
    
    \item $J = \left\{ \Lambda^k \left[ n \right] \rightarrowtail \Delta
    \left[ n \right] | 0 \leq k \leq n, n \in \mathbb{N} \setminus \left\{ 0
    \right\} \right\}$.
  \end{itemizeminus}
\end{theorem}

\begin{definition}
  Let $\C$ be a category with coproducts. A {\tmem{cylinder}} $\is = ( I,
  \iota_0, \iota_1, \sigma ) $on $\C$ is given by the following
  data:{\tmstrong{}}
  \begin{itemizeminus}
    \item an endofunctor $I : \C \longrightarrow \C$;
    
    \item natural transformations $\iota^0, \iota^1 : \tmop{id}_{\C}
    \Rightarrow I$ and $\sigma : I \Rightarrow \tmop{id}_{\C}$ such that
    $\sigma \circ i^0 = \sigma \circ \iota^1 = \tmop{id}_{\tmop{id}_{\C}}$ ;
  \end{itemizeminus}
  A cylinder is {\tmem{cartesian}} if
  \begin{enumerateroman}
    \item $I$ preserves monos;
    
    \item the canonical morphism $[ \iota^0_C, \iota^1_C ] : C + C
    \longrightarrow I ( C )$ is mono for all $C \in \C$;
    
    \item the naturality square
    
    \begin{center}
      $
\xy
\xymatrix{ 
K \ar@{ >->}[d]_j \ar[r]^{\iota^w_K} &
I(K) \ar@{ >->}[d]^{I(j)}
\\
L \ar[r]_{\iota^w_L} & I(L)
}
\POS( 5 , -4.5);\POS( 5,-2.2 ) \connect@{-}
\POS( 5 , -4.5);\POS( 2.7,-4.5 ) \connect@{-}
\endxy
$

    \end{center}

    {\noindent}is a pullback square for all monos $j$ and $w \in \{ 0, 1 \}$.
  \end{enumerateroman}
\end{definition}

\begin{definition}
  \label{def-homdat}Let $\C$ be a small category and $\widehat{\C}  \deq 
  \sett^{\C^{\tmop{op}}}$ its category of presheaves. An {\tmem{elementary
  homotopical datum}} on $\C$ is a cartesian cylinder $\is = ( I, \iota_0,
  \iota_1, \sigma )$ on $\widehat{\C}$ such that $I$ preserves colimits. A
  {\tmem{homotopical datum}} on $\C$ is a pair $( \is, S )$ with $\is$ a
  homotopical datum on $\C$ and $S \subseteq \widehat{\C}_1$ a set of monos. 
\end{definition}

As the name suggests, an elementary homotopical datum gives a notion of
homotopy on morphisms of presheaves.

\begin{proposition}
  Let $\C$ be a small category and $\is$ an {\tmem{elementary homotopical
  datum}} on $\C$. Given morphisms of presheaves $f_0, f_1 : X \longrightarrow
  Y$ let
  \[ f_0 \sim^1 f_1  \dee \exists h : I ( X ) \longrightarrow Y. h \circ
     \iota^w_X = f_w \]
  for $w \in \{ 0, 1 \}$. The equivalence relation $\sim^{\is}$on
  $\widehat{\C_1}$ generated by $\sim^1$is a congruence.
\end{proposition}

\begin{definition}
  Let $\C$ be a small category, $\is$ an {\tmem{elementary homotopical datum}}
  on $\C$ and $j : K \rightarrowtail L$ a mono in $\widehat{\C}$.
  \begin{enumerateroman}
    \item $\Theta ( j )$ is the comparison morphism in

    \begin{center}
      $
\xy
\xymatrix{
K \ar@{ >->}[d]_j \ar@{ >->}[r]^{\delta^w_K} 
&
I(K) \ar[d] \ar@/^1pc/@{ >->}[ddr]^{I(j)}
&
\\
L \ar[r]  \ar@/_1pc/@{ >->}[drr]_{\delta^w_L} & 
\bullet \ar@{.>}[dr]_{\Theta(j)}
&
\\
&& I(L) 
}

\POS( 11.7 , -9);\POS( 11.7,-11.3 ) \connect@{-}
\POS( 11.7 , -9);\POS( 14,-9 ) \connect@{-}
\endxy
$

    \end{center}
    
    \item $\Lambda ( j )$ is the comparison morphism from in

    \begin{center}
      $
\xy
\xymatrix {
K + K \ar@{ >->}[d]_{j+j} 
\ar@{ >->}[r]^-{ [\delta^0_K,\delta^1_K]} 
&
I(K) \ar[d] \ar@/^1pc/@{ >->}[ddr]^{I(j)}
&
\\
L+L \ar[r]  
\ar@/_1pc/@{ >->}[drr]_{ [\delta^0_L,\delta^1_L]} & 
\bullet \ar@{.>}[dr]_{\Lambda(j)}
&
\\
&& I(L) 
}

\POS( 14 , -9);\POS( 14,-11.3 ) \connect@{-}
\POS( 14 , -9);\POS( 16.3,-9 ) \connect@{-}
\endxy
$

    \end{center}
  \end{enumerateroman}
  Given a set of monos $M \in \widehat{\C}$, let $\Theta ( M ) \deq \{ \Theta
  ( j ) | j \in M \}$ and $\Lambda ( T ) \deq \{ \Lambda ( j ) | j \in M \}$.
\end{definition}

\begin{theorem}
  \label{theo-cisi} {\dueto{Cisinski}}Let $\C$ be a small category, $\text{$(
  \is, S )$}$ be a homotopical datum on $\C$ and $\ms \in \widehat{\C_1}$ be a
  set of monos such that $\llp{\rlp{\ms}}$ is the class of all monos. Let
  \begin{itemizeminus}
    \item $\Lambda_0  \deq S \cup \Theta ( \ms )$ and $\Lambda_{n + 1}  \deq
    \Lambda ( \Lambda_n )$;
    
    \item $\Lambda_{\is} ( S, \ms ) \deq  \bigcup_{n \geqslant 0} \Lambda_n$.
  \end{itemizeminus}
  $\widehat{\C}$ admits a cofibrantly generated model structure where the
  cofibrations are the monos and the weak equivalences are the morphisms $f :
  X \longrightarrow Y$ inducing a bijection
  \[ f^{\ast} : \left(  \widehat{\C} / \sim^{\is} \right) ( Y, T )
     \cong  \left(  \widehat{\C} / \sim^{\is} \right) ( X, T )
  \]
  for all $T \in \C$ such that $T \longrightarrowlim^{!_T} 1 \in
  \rlp{\Lambda_{\is} ( S, \ms )}$. 
\end{theorem}

Theorem \ref{theo-cisi} works in fact for all topoi, not only those of
presheaves {\cite{cisinski-topos}}. \\

\begin{theproofof}
\textbf{theorem \ref{theo-quillen-sset}}. Set
\begin{itemizeminus}
  \item $S = \varnothing$,
  
  \item $\ms  \deq  \left\{ \partial \left[ n \right] \rightarrowtail \Delta
  \left[ n \right] | n \in \mathbb{N} \right\}$ and
  
  \item $\is  \deq (_- ) \times \Delta [ 1 ]$
\end{itemizeminus}
and apply theorem \ref{theo-cisi}. 
\end{theproofof}

However, as far as labor is concerned, there is no 
thing like a free lunch.
What one spares with the existence is spent with the characterisations of
$\ws$ and $J$ ($I$ is easy). Nonetheless, \ref{theo-quillen-sset} 
theorem allows to isolate the
non-structural part of a task at hand.

\subsection{Locally Presentable Categories for the Homotopy Theorist}

\begin{definition}
  Suppose $\mathbb{A}$ has all coproducts. A family of objects $\left( G_i
  \right)_{i \in I}$ is a family of generators if the comparison morphism
  \[ \gamma_C \deq \left[ f \right]_{i \in I, \hspace{0.25em} f \in \mathbb{A}
     \left( G_i, C \right)} : \hspace{0.5em} \left( \coprod_{i \in I,
     \hspace{0.25em} f \in \mathbb{A} \left( G_i, C \right)} G_i \right)
     \rightarrow C \]
  is epi for all $C \in \mathbb{A}$. A family of generators is
  \begin{enumerateroman}
    \item {\tmem{strong}} if $\gamma_C \in \llp{\tmop{Monos}}$ for all $C \in
    \mathbb{A}$;
    
    \item {\tmem{dense}} if, given the full subcategory $\mathbbm{G} \subseteq
    \mathbbm{A}$ such that $\mathbbm{G}_0 = \left( G_i \right)_{i \in I}$,
    $\left( C, \hspace{0.25em} \left( f \right)_{f \in \mathbb{G} / C}
    \right)$ is a colimit of $\tmop{dom} : \hspace{0.25em} \mathbb{G} / C
    \rightarrow \mathbb{A}$ for all $C \in \mathbb{A}$.
  \end{enumerateroman}
  A one-member family of generators is called a generator (respectively a
  strong generator, respectively a dense generator).
\end{definition}

A familiar example is given by the Yoneda embedding: the family of all
representable presheaves $\left( \mathbb{B} \left( \_, B \right) \right)_{B
\in \mathbb{B}_0}$ over some category $\mathbb{B}$ is a dense generating
family in $\sett^{\mathbb{B}^{op}}$.

\begin{definition}
  Let $\alpha$ be a regular cardinal. $C \in \mathbb{A}$ is
  $\alpha$-\textit{presentable} provided $\mathbb{A} \left( C, \_ \right)$
  preserves $\alpha$ \textit{-}filtered colimits. It is \textit{presentable}
  if there is an $\alpha$ sucht that it is $\alpha$-presentable. 
\end{definition}

An $\alpha$-presentable $C \in \mathbb{A}$ is $\beta$-presentable for any
regular $\beta < \alpha$. Finitely presentable groups are presentable.
Presentable topological spaces are precisely the discrete ones i.e. there is
no regular cardinal $\alpha$ for which a topological space is
$\alpha$-presentable. Gabriel and Ulmer observe that ``...the presentable
individuals are the discrete ones, an exemplary society!''
\cite[p.64]{GabrielP:locpc} {\footnote{``Insbesondere sind die
pr\"asentierbare Individuen gereade die Diskreten, eine vorbildliche
Gesellschaft!''}}.

\begin{definition}
  Let $\alpha$ be a regular cardinal. The category $\mathbb{A}$ is
  \textit{locally $\alpha$-presentable} provided
  \begin{enumerate}
    \item $\mathbb{A}$ is cocomplete;
    
    \item $\mathbb{A}$ has a strong family of generators $\left( G_i
    \right)_{i \in I}$;
    
    \item each $G_i$ is $\alpha$-presentable. 
  \end{enumerate}
\end{definition}

\begin{theremark}
  \label{rem-2catlocpres} {\tcat} is locally presentable. It is cocomplete by
  proposition \ref{pro:2cat-complete} and it is easy to see that the
  2-category $\tmmathbf{W}_2$ given by

  \begin{center}
    $
\xymatrix{
X \rtwocell^f_g{\alpha} & Y
}
$

  \end{center}

  {\noindent}(a.k.a ``the walking 2-cell'' or ``the free-living 2-cell'') is a
  strong $\aleph_0$-presentable generator.
\end{theremark}

\begin{proposition}
  \label{pro:lpc1}Let $\alpha$ be a regular cardinal and $\mathbb{A}$ be
  locally $\alpha$-presentable. Let $\mathbb{G}$ be the full subcategory
  spanned by $\mathbb{A}$'s generating family $\left( G_i \right)_{i \in I}$. 
Then
  \begin{enumerate}
    \item The closure $\mathbb{P}$ of $\mathbb{G}$ under $\alpha$-colimits
    exists and is equivalent to a small category;
    
    \item $\mathbb{P}$'s $\alpha$-colimits are computed as in $\mathbb{A}$;
    
    \item every object in $\mathbb{P}$ is $\alpha$-presentable;
    
    \item $\mathbb{P}_0$ is a dense generator in $\mathbb{A}$.
  \end{enumerate}
\end{proposition}

\begin{proposition}
  Let $\alpha$ be a regular cardinal and $\mathbb{A}$ be locally
  $\alpha$-presentable. For every $C \in \mathbb{A}$ there is a regular
  cardinal $\alpha_C$ such that $C$ is $\alpha_C$-presentable.
\end{proposition}

\begin{corollary}
  \label{coro-small}The small object argument applies to any set $I \subseteq
  \mathbb{A}_1$.
\end{corollary}

Corollary \ref{coro-small} is the main reason for the interest of homotopy
theorists in locally presentable categories.

\subsection{Creation of Model Structures by Right Adjoints}

\begin{definition}
  \label{def:creation}Let $\mathbb{M}$ be a model category and
  
  \begin{center}
    $\xymatrix{ 
\C \rtwocell^{F}_{U}{`{\perp}} & \M
}$

  \end{center}
  
  {\noindent}be an adjunction. $U$ \textit{creates} a model structure on
  $\mathbb{C}$ if there is a model structure on $\mathbb{C}$ such that
  $\mathcal{F}_{\mathbb{C}} = U^{- 1} \left( \mathcal{F}_{\mathbb{M}} \right)$
  and $\mathcal{W}_{\mathbb{C}} = U^{- 1} \left( \mathcal{W}_{\mathbb{M}}
  \right)$.
\end{definition}

\begin{proposition}
  \label{prop-creation}Let $\M$ be a cofibrantly generated model category with
  $I$ and $J$ the sets of generating cofibrations and acyclic cofibrations,
  respectively. Let $F \dashv U$ and $\C$ be as in definition
  \ref{def:creation}. Suppose
  \begin{enumerateroman}
    \item $\dom  \left( F ( i ) \right)$ is small with respect to $F ( I )$
    for all $i \in I$ and $\dom  \left( F ( j ) \right)$ is small with respect
    to $F ( J )$ for all $j \in J$;
    
    \item the composition of any $( \lambda, \ws_{\M} )$-sequence $\lambda
    \longrightarrow \M$ is a weak equivalence for all $\lambda \in
    \tmmathbf{\ord}$;
    
    \item $U$ preserves colimits of $\lambda$-sequences for all $\lambda \in
    \tmmathbf{\ord}$; and
    
    \item for every $A \longrightarrowlim^j B \in J$ and for every pushout

    \begin{center}
      $
\xy
$\xymatrix{ 
F(A)\ar [d]_{F(j)} \ar [r] ^f & C\ar [d]^g\\
F(B)\ar [r] & F(B) +_{F(A)} C
}$

\POS( 18.5 , -7.5);\POS( 18.5,-9.9 ) \connect@{-}
\POS( 18.5 , -7.5);\POS( 20.9,-7.5 ) \connect@{-}
\endxy
$

    \end{center}
\mbox{}\\\\
    {\noindent}in $\C$, the morphism $U ( g )$ is a weak equivalence.
  \end{enumerateroman}
  Then the adjoint pair $F \dashv U$ creates a cofibrantly generated model
  category structure on $\C$, where $F ( I )$ and $F ( J )$ are the generating
  cofibrations and generating acyclic cofibrations, respectively.
\end{proposition}

Proposition \ref{prop-creation} is an easy consequence of Kan's Theorem on
creation of model category structure \cite[thm. 11.3.2]{hirsch}.

\section{\label{sec:modelstr}A Model Structure \`a la Thomason}

This section essentially revisits and generalizes categorical techniques
developed by Fritsch and Latch {\cite{frilat}}, Thomason {\cite{thoma-model}}
and Cisinski {\cite{cisinski-thoma,cisinski-thesis}}. However, it turns out
that not everything carries over by tagging a ``2-'' in front. It is crucially
the case for Cisinski's ``immersions'', a generalization of Thomason's
``Dwyer-morphisms''. We call the relevant 2-categorical notion ``skew
immersion''.

\subsection{Ideals in Categories}

\begin{definition}
  Let $\mathcal{\A} \subseteq \mathcal{\B}$ be an inclusion of categories.
  $\mathcal{\A}$ is an {\tmem{L-ideal}} in $\B$ if
  \[ \forall f \in \B_1 . \cod \left( f \right) \in \A_0 \hspace{0.25em}
     \Rightarrow f \in \A_1 \]
  and $\mathcal{}$an \textit{R-ideal} in $\B$ if
  \[ \forall f \in \B_1 . \dom \left( f \right) \hspace{0.25em} \in \A_0
     \Rightarrow f \in \A_1 \]
\end{definition}

In the literature,  $L$-Ideals are called {\tmem{left ideals}},
{\tmem{sieves}} or {\tmem{cribles}} while $R$-Ideals are called {\tmem{right
ideals}}, {\tmem{cosieves}} or {\tmem{cocribles}}
{\cite{frilat,thoma-model,cisinski-thoma,cisinski-thesis}}

\begin{definition}
  Let $\I$ be the category generated by $L \longrightarrowlim^t R$ and
  $\iota^L, \iota^R : 1 \longrightarrow \I$ be the global elements of $\I$
  with image generated by $L$ respectively by $R$. Let further $\partial^L 
  \deq  \cod$, $\partial^R  \deq  \dom$ and
  \[ \begin{array}{llll}
       \overline{(_- )} : & \{ L, R \} & \longrightarrow & \{ L, R \}\\
       & L & \longmapsto & R\\
       & R & \longmapsto & R
     \end{array} \]
  be the toggling map.
\end{definition}

\begin{proposition}
  \label{prop-catideals}Let $\mathcal{\A} \subseteq \mathcal{\B}$ be an
  inclusion of categories and $\nu \in \{ L, R \}$. The following are
  equivalent.
  \begin{enumerateroman}
    \item $\A$ is a $\nu$-ideal;
    
    \item there is a functor $\chi_{\A} : \B \longrightarrow I$ such that $\A
    \cong \chi_{\A}^{\ast} ( \iota^{\nu} )$;
    
    \item $\mathcal{\A} \subseteq \mathcal{\B}$ is a full inclusion and there
    is a commuting square 
\[\xymatrix{
\A \ar@{ >->}[d] \ar[r]^{!} & 1 \ar@{ >->}[d]^{\iota^\nu}
\\
\B \ar[r]_{\chi_\A} & \I
}\]
    {\noindent}such that
    \[ ( \chi_{\A} ) ( B ) = \left\{\begin{array}{ll}
         \nu & B \in \A_0\\
         \bar{\nu} & B \in \B_0 \backslash \A_0
       \end{array}\right.  \]
  \end{enumerateroman}
\end{proposition}

\begin{theproof}
  ({\tmem{i}})$\Rightarrow$({\tmem{ii}})  The functor given by
  \[ \begin{array}{llll}
       ( \chi_{\A} )_{_0} : & B & \longmapsto & \left\{\begin{array}{ll}
         \nu & B \in \A_0\\
         \bar{\nu} & B \in \B_0 \backslash \A_0
       \end{array}\right.\\
       &  &  & 
     \end{array} \]
  and
  \[ \begin{array}{llll}
       ( \chi_{\A} )_{_1} : & f & \longmapsto & \left\{\begin{array}{ll}
         \id_{\nu} & f \in \A_1 \\
         \theta & \partial^{\bar{\nu}} ( f ) \in \A_0 \wedge \partial^{\nu} (
         f ) \in \B_0 \backslash \A_0\\
         \id_{\bar{\nu}} & f \in \B_1 \backslash \A_1
       \end{array}\right.
     \end{array} \]
  is well-defined since $\A$ is a $\nu$-ideal. It is immediate that $\A \cong
  \chi_{\A}^{\ast} ( \iota^{\nu} )$.
  
  ({\tmem{ii}})$\Rightarrow$({\tmem{iii}}) Let $f \in \B_1$ such that
  $\partial^{\nu} ( f ) \in \chi_{\A}^{\ast} ( \iota^{\nu} )_{_0}$ and
  $\partial^{\bar{\nu}} ( f ) \in \chi_{\A}^{\ast} ( \iota^{\nu} )_{_0}$. Then
  $f \in \chi_{\A}^{\ast} ( \iota^{\nu} )_{_1}$ by the underlying graph
  structure, so the inclusion is full.
  
  ({\tmem{iii}})$\Rightarrow$({\tmem{i}}) Let $f \in \B_1$ such that
  $\partial^{\nu} ( f ) \in \A_0$. Then $\chi_{\A} ( \partial^{\nu} ( f ) ) =
  \nu$ by definition of $\chi_{\A} $and $\chi_{\A} ( \partial^{\bar{\nu}} ( f
  ) ) = \nu$ by the underlying graph structure, hence  $\partial^{\bar{\nu}} (
  f ) \in \A_0$. But $\A$ is a full subcategory so $f \in \A_1$.
\end{theproof}

\begin{definition}
  The functor $\chi_{\A}$ of proposition \ref{prop-catideals} is called the
  ideal's {\tmem{characteristic morphism}}.
\end{definition}

\begin{theremark}
  An ideal is in particular always a full subcategory. A characteristic
  morphism in necessarily unique.
\end{theremark}

\subsection{Ideals in 2-Categories}

The notion of ideal carries over as expected to 2-categories.

\begin{definition}
  Let $\mathcal{A} \subseteq \mathcal{B}$ be an inclusion of 2-categories.
  $\mathcal{A}$ is an \textit{{\tmem{L-ideal}}} in $\bs$ if
  \[ \forall \alpha \in \mathcal{B}_2 . ( \cod \circ \dom ) \left( \alpha
     \right) = \left( \cod \circ \cod \right) \left( \alpha \right) \in \as_0
     \hspace{0.25em} \Rightarrow \alpha \in \mathcal{A}_2 \]
  and $\mathcal{}$an \textit{R-ideal} in $\bs$ if 
  \[ \forall \alpha \in \mathcal{B}_2 . ( \dom \circ \dom ) \left( \alpha
     \right) = \left( \dom \circ \cod \right) \left( \alpha \right) \in \as_0
     \hspace{0.25em} \Rightarrow \alpha \in \mathcal{A}_2 \]
\end{definition}

We also call $L$-ideals {\tmem{left ideals}} and $R$-ideals {\tmem{right
ideals}}.

\begin{proposition}
  \label{prop-2catideals}Let $\is$ be the 2-category with trivial 2-cells such
  that $\left\lfloor \is \right\rfloor = \I$. Let $\nu \in \{ L, R \}$ and
  $\mathcal{A} \subseteq \mathcal{B}$ be an inclusion of 2-categories  The
  following are equivalent.
  \begin{enumerateroman}
    \item $\as$ is a $\nu$-ideal in $\bs$;
    
    \item $\left\lfloor \as \right\rfloor$ is a $\nu$-ideal in $\left\lfloor
    \bs \right\rfloor$ and $\text{$\mathcal{A} \subseteq \mathcal{B}$}$ is a
    locally full inclusion;
    
    \item there is a 2-functor $\text{$\chi_{\as} : \bs \longrightarrow I$}$
    such that $\left\lfloor \chi_{\as} \right\rfloor = \chi_{\left\lfloor \as
    \right\rfloor}$.
  \end{enumerateroman}
\end{proposition}

\begin{theproof}
  ({\tmem{i}})$\Rightarrow$({\tmem{ii}}) $\left\lfloor \as \right\rfloor$ is a
  $\nu$-ideal by instantiating the defintion on the identity 2-cells. Suppose
  $\dom_1 ( \alpha ) \in \as_1$ and $\cod_1 ( \alpha ) \in \as_1$. Then in
  particular $( \partial^{\nu} \circ \dom_1^{} ) ( \alpha ) \in \as_0$ and  $(
  \partial^{\nu} \circ \cod_1^{} ) ( \alpha ) \in \as_0$, hence $\alpha \in
  \as_2$.
  
  ({\tmem{ii}})$\Rightarrow$({\tmem{i}}) Suppose $( \partial^{\nu} \circ
  \dom_1^{} ) ( \alpha ) \in \as_0$ and  $( \partial^{\nu} \circ \cod_1^{} ) (
  \alpha ) \in \as_0$. Then $\dom_1 ( \alpha ) \in \as_1$ and $\cod_1 ( \alpha
  ) \in \as_1$ since $\left\lfloor \as \right\rfloor$ is a $\nu$-ideal. But
  $\text{$\mathcal{A}$}$ is a locally full sub2-category so $\alpha \in
  \as_2$.
  
  ({\tmem{ii}})$\Leftrightarrow$({\tmem{iii}}) Obvious.
\end{theproof}

\begin{definition}
  The 2-functor $\text{$\chi_{\as} : \bs \longrightarrow \is$}$ of proposition
  \ref{prop-2catideals} is called the ideal's {\tmem{characteristic
  morphism}}.
\end{definition}

\begin{theremark}
  An ideal inclusion is in particular always full and locally full. The
  characteristic morphism is necessarily unique and $\as \cong
  \chi_{\as}^{\ast} ( i^{\nu} )$.
\end{theremark}

\begin{lemma}
  \label{lem-idealstab}Ideals are stable under pullback and pushout.
\end{lemma}

\begin{theproof}
  Let $\nu \in \{ L, R \}$. The first assertion follows immediately from the
  pullback lemma:

  \begin{center}
    $
\xy
\xymatrix{
\as' \ar[r] \ar@{ >->}[d] & \as \ar[r] \ar@{ >->}[d] & 
1 \ar@{ >->}[d]^{\iota^\nu}
\\
\bs' \ar[r] & \bs \ar[r]_{\chi_\as} & \mathcal{I} 
}
\POS( 4.5 , -4);\POS( 4.5,-1.7 ) \connect@{-}
\POS( 4.5 , -4);\POS( 2.2,-4 ) \connect@{-}
\POS( 18.5 , -4);\POS( 18.5,-1.7 ) \connect@{-}
\POS( 18.5 , -4);\POS( 16.2,-4 ) \connect@{-}
\endxy
$

  \end{center}

  {\noindent}For the second, consider the diagram

  \begin{center}
    $
\xy
\xymatrix{
\as \ar[r] \ar@{ >->}[d] &
\as' \ar@{ >->}[d] \ar[dr]^{!} & 
\\
\bs \ar[r] \ar@/_/[drr]_{\chi_\bs} &
\bs' \ar@{.>}[dr]^<<<<{\chi_{\bs'}} & 1 \ar@{ >->}[d]^{\iota^\nu}
\\
&& \mathcal{I}
}

\POS( 10 , -9.4);\POS( 10,-11.4 ) \connect@{-}
\POS( 10 , -9.4);\POS( 12,-9.4 ) \connect@{-}
\endxy
$

  \end{center}

  {\noindent}with $\chi_{\bs'}$ given by universal property. By remark
  \ref{rem-pushinc}, $\as' \subseteq \bs'$ is full and locally full and the
  pushout square is

  \begin{center}
    $
\xy
\xymatrix @*[r]{
\as_0 \ar[rr]^{f_0} \ar@{ >->}[d] & &
\as_0' \ar@{ >->}[d]^{\text{in}_{\as_0'}}
\\
\bs_0 \ar[rr]_{f_0 + \text{id}_{(\mathcal{B}_0 \setminus \mathcal{A}_0})} &&
\as_0' +_{\as_0} \left(\bs_0 \setminus \as_0\right)
}
\POS( 31.5 , -9.5);\POS( 31.5,-11.8 ) \connect@{-}
\POS( 31.5 , -9.5);\POS( 33.8,-9.5 ) \connect@{-}
\endxy
$

  \end{center}

  {\noindent}on objects. We have
  \[ \chi_{\mathcal{B'}} |_{\as_0'} = \iota^{\nu} \circ !_{\as'_0} \]
  and
  \[ ( \chi_{\bs'} \circ ( f_0 + \tmop{id}_{( \bs_0 \backslash \as_0 )} ) )
     |_{\bs_0 \backslash \as_0} = \chi_{\bs'} |_{\bs_0 \backslash \as_0} =
     \chi_{\bs} |_{\bs_0 \backslash \as_0} \]
  hence
  \[ \chi_{\bs'} ( B' ) = \left\{\begin{array}{ll}
       \nu & B' \in \as'_0\\
       \bar{\nu} & B' \in ( \bs'_0 \backslash \as'_0 ) \cong ( \bs_0
       \backslash \as_0 )
     \end{array}\right. \]
  so the assertion follows by proposition \ref{prop-catideals}.
\end{theproof}

\begin{definition}
  Let $\as$ be a 2-category and $X \subseteq \as_0$. $\left. \left.
  \right\lceil X \right\rceil \subseteq \as$ is the full and locally full
  sub-2-category such that $\left\lceil X \right\rceil_0 = X$.
\end{definition}

\begin{lemma}
  \label{lem:vwb}Let $\as \subseteq \ws \subseteq \bs$ be inclusions of
  2-categories with $\as$ a left ideal and $\ws$ a right ideal. Let $\V  \deq 
  \left\lceil \bs_0 \setminus \as_0 \right\rceil$. The image of the pullback
  square

  \begin{center}
    $
\xy
\xymatrix{
(\V) \cap \ws \ar@{ >->}[r] \ar@{ >->}[d] & \ws \ar@{ >->}[d] \\
\V \ar@{ >->}[r] & \bs
}
\POS( 6.1 , -6.3);\POS( 6.1,-3.9 ) \connect@{-}
\POS( 6.1 , -6.3);\POS( 3.7,-6.3 ) \connect@{-}
\endxy
$

  \end{center}
  
  {\noindent}under $N_2$ is a pushout square.
\end{lemma}

\begin{theproof}
  The comparison map $c$ is an injection for all $n \in \N$:

  \begin{center}
    $
\xy
\xymatrix{
N_2((\V) \cap \ws)_n \ar@{ >->}[d] \ar@{ >->}[r] 
&
N_2(\ws)_n \ar@{ >->}[d] \ar@/^1pc/@{ >->}[ddr]
&
\\
N_2(\V)_n \ar@{ >->}[r]  \ar@/_1pc/@{ >->}[drr] & 
\bullet \ar@{ >.>}[dr]_{c_n}
&
\\
&& N_2(\bs)_n
}
\POS( 25 , -9.5);\POS( 25,-11.9 ) \connect@{-}
\POS( 25 , -9.5);\POS( 27.4,-9.5 ) \connect@{-}
\endxy
$

  \end{center}

  {\noindent}Recall that $N_2 ( \as )_n = \tmmathbf{\nlax} ( [ n ], \as )$
  (c.f. remark \ref{rem:2-nerve}). Suppose $n = 0$. We have
  \[ N_2 ( \V )_0 +_{N_2 ( ( \V ) \cap \ws )} N_2 ( \ws )_0 \cong ( \V )_0
     +_{( ( \V )_0 \cap \ws_0 )}  \ws_0 \cong ( \bs_0 \backslash \as_0 ) \cup
     \ws_0 = \bs_0 \]
  since $\as_0 \subseteq \ws_0$. In particular, $c_0$ is a surjection. Suppose
  now $n > 0$. $\V$ is a right ideal since $\as$ is a left ideal and $\ws$ is
  a right ideal by hypothesis, hence the image of  a lax functor $[ n ]
  \longrightarrow \bs$ is in $\V$ or in $\ws$ so $c_n$ is a surjection for all
  $n \in \N$.
\end{theproof}

\begin{definition}
  \label{def-thoquot}Let $\as \rightarrowtail \bs $ be an inclusion of
  2-categories. The 2-category $\bs / \as$ is given by the pushout square

  \begin{center}
    $
\xy
\xymatrix{
\as \ar[r]^{!_\as} \ar@{ >->}[d] & 1 \ar@{ >->}^{\rho_\as}[d]
\\
\bs \ar[r]_{\rho_\bs} & \bs / \as
}
\POS( 10 , -8);\POS( 10,-10.4 ) \connect@{-}
\POS( 10 , -8);\POS( 12.4,-8 ) \connect@{-}
\endxy
$

  \end{center}

\end{definition}

\begin{proposition}
  \label{prop-thoquot}Let $I : \as \rightarrowtail \bs $ be an inclusion of
  2-categories, $F : \as \longrightarrow \as'$ a 2-functor and

  \begin{flushright}
    \begin{center}
      $
\xy
\xymatrix{
\as \ar[r]^F \ar@{ >->}[d] & \as' \ar@{ >->}[d]
\\
\bs \ar[r] & \as' +_\as \bs
}
\POS( 13.5 , -8);\POS( 13.5,-10.3 ) \connect@{-}
\POS( 13.5 , -8);\POS( 15.8,-8 ) \connect@{-}
\endxy
$

    \end{center}
  \end{flushright}

  {\noindent}the corresponding pushout square. Then
  \[ \bs / \as \cong \left( \as' +_{\as}  \bs \right) / \as' \]
\end{proposition}

\begin{theproof}

  \begin{center}
    $
\xy
\xymatrix @*[r]{
\as \ar[r]^F \ar@{ >->}[d] 
&
\as' \ar@{ >->}[d] \ar[r]^\!
&
1 \ar@{ >->}[d]
\\
\bs \ar[r] 
& 
\as' +_\as \bs \ar[r]
&
\bs / \as \cong (\as' +_\as \bs) / \as'
}
\POS( 13.5 , -9);\POS( 13.5,-11.4 ) \connect@{-}
\POS( 13.5 , -9);\POS( 15.9,-9 ) \connect@{-}
\POS( 47 , -9);\POS( 47,-11.4 ) \connect@{-}
\POS( 47 , -9);\POS( 49.4,-9 ) \connect@{-}
\endxy
$

  \end{center}

\end{theproof}

\begin{corollary}
  \label{cor-thoquot} $\bs \backslash \as \cong \left( \as' +_{\as}  \bs
  \right) \backslash \as'$ provided $\as$ is an ideal.
\end{corollary}

\subsection{Distorsions}

\begin{definition}
  Let $\kappa^{\nu} : \as \cong \as \times 1 \longrightarrowlim^{\id \times
  \iota^{\nu}} \as \times \is$ for $\nu \in \{ L, R \}$. Let $F, G :
  \hspace{0.25em} \as \rightarrow \bs$ be 2-functors. A {\tmem{distortion}}
  \textit{$\varepsilon : \hspace{0.25em} F \rightsquigarrow G$} is given by a
  normal lax functor $\underline{\varepsilon} : \hspace{0.25em} \as \times \is
  \longrightarrow \bs$ such that

  \begin{center}
    $\xymatrix{
\as \ar[dr]_F \ar[r]^-{\kappa^l} & \as \times \mathcal{I}
\ar[d]^{\underline{\varepsilon}} & \as \ar[l]_-{\kappa^r} 
\ar [dl]^G \\
& \bs & 
}$

  \end{center}
  
  {\noindent}commutes in $\tcatln$.
\end{definition}

\begin{theremark}
  \label{rem-simho} $N_2$ extends to a product-preserving functor
  \[ \widetilde{N_2} : \hspace{0.25em} \widetilde{\tcat} \rightarrow \sset \]
  It follows that a distortion $\varepsilon : \hspace{0.25em} F
  \rightsquigarrow G$ gives rise to a simplicial homotopy $N_2 \left( F
  \right) \simeq N_2 \left( G \right)$.
\end{theremark}

\begin{proposition}
  Let $F, G : \hspace{0.25em} \as \rightarrow \bs$ be 2-functors. The
  following are equivalent.
  \begin{enumerateroman}
    \item There is a distortion $\varepsilon : \hspace{0.25em} F
    \rightsquigarrow G$;
    
    \item there are
    \begin{itemizeminus}
      \item a morphism $\varepsilon_f : \hspace{0.25em} F ( A )
      \longrightarrow G ( B )$ for all $f \in \as_1$;
      
      \item a 2-cell $\varepsilon_{\alpha} : \hspace{0.25em} \varepsilon_f
      \Longrightarrow \varepsilon_{f'}$ for all 
$\alpha: f \Rightarrow f' \in \as_2$;
      
      \item 2-cells $\text{$\varepsilon_{f, g}^L : \hspace{0.25em}
      \varepsilon_g \circ F ( f ) \Rightarrow \varepsilon_{g \circ f}$}$ and
      
      \item $\varepsilon_{f, g}^R : \hspace{0.25em} G ( g ) \circ
      \varepsilon_f \Rightarrow \varepsilon_{g \circ f}$  for all composable
      $f, g \in \as_1$,
    \end{itemizeminus}
    such that
    \begin{description}
      \item[lf1] $\varepsilon_{( \beta \bullet \alpha )} = \varepsilon_{\beta}
      \bullet \varepsilon_{\alpha}$ for all vertically composable  $\alpha,
      \beta \in \as_2$;
      
      \item[lf2] $\varepsilon_{\id_f} = \id_{\varepsilon_f}$ for all $f \in
      \as_1$;
      
      \item[n1] $\varepsilon^R_{\cod^1 ( \varphi ), \cod^1 ( \theta )} \bullet
      \left( G ( \theta ) \circ \varepsilon_{\varphi} \right) =
      \varepsilon_{\theta \circ \varphi} \bullet \varepsilon^R_{\dom^1 (
      \varphi ), \dom^1 ( \theta )}$ and
      
      \item[n2] $\varepsilon^L_{\cod^1 ( \varphi ), \cod^1 ( \theta )} \bullet
      \left( \varepsilon_{\theta} \circ F ( \varphi ) \right) =
      \varepsilon_{\theta \circ \varphi} \bullet \varepsilon^L_{\dom^1 (
      \varphi ), \dom^1 ( \theta )}$ for all horizontally composable $\varphi,
      \theta \in \as_2$;
      
      \item[c1] $\left. \varepsilon^R_{g \circ f, h} \bullet \left( G ( h
      \right) \circ \varepsilon^R_{f, g} \right) = \varepsilon^R_{f, g \circ
      h}$,
      
      \item[c2] $\varepsilon^L_{f, h \circ g} \bullet \left( \varepsilon^L_{g,
      h} \circ F ( f ) \right) = \varepsilon^L_{g \circ f, h}$ and
      
      \item[c3] $\varepsilon^R_{g \circ f} \bullet \left( G ( h ) \circ
      \varepsilon^L_{f, g} \right) = \varepsilon^L_{f, h \circ g} \bullet
      \left( \varepsilon^R_{g, h} \circ F ( f ) \right)$ for all composable
      $f, g, h \in \as_1$.
    \end{description}
  \end{enumerateroman}
\end{proposition}

\begin{theproof}  Let $A \in \as_0$, $f \in \as_1$ and $\nu \in \{ L, R \}$.
The values of $\underline{\varepsilon}$ on $( A, \nu )$ respectively on $( f,
\id_{\nu} )$ are determined by $F$ for $\nu = L$ and by $G$ for $\nu = R$. The
associated structural 2-cells are all trivial. On the other hand,
\[ \varepsilon_f  \deq  \underline{\varepsilon} ( f, t ) \]
and
\[ \varepsilon_{\alpha}  \deq  \underline{\varepsilon} ( \alpha, \id_t ) \]
are the remaining values while
\[ \varepsilon^L_{f, g}  \deq \gamma_{( f, \id_L ), ( g, t )} \]
and
\[ \varepsilon^R_{f, g}  \deq \gamma_{( f, t ), ( g, \id_R )} \]
are the remaining structural 2-cells.
\end{theproof}

\begin{theremark}
  \label{rem:wk}Distortions do not compose in general (neither vertically nor
  horizontally), yet they can be whiskered on the left as well as on the
  right. 
\end{theremark}

\begin{theremark}
  \label{rem-dist-conv}Some instances of the equations governing a distortion
  become conveniently simpler. Let $\varepsilon : F \rightsquigarrow G$ be a
  distortion. Given

  \begin{center}
    $
\xymatrix{
A \ar[r]^{\id_A} & A \rtwocell^u_v{\theta} & B \ar[r]^{\id_B} & B
}
$

  \end{center}

  {\noindent}let  $\varepsilon_A  \deq \varepsilon_{\id_A}$,
  $\varepsilon^{\nu}_{A, u}  \deq \varepsilon^{\nu}_{\id_A, u}$ and
  $\varepsilon^{\nu}_{u, B}  \deq \varepsilon^{\nu}_{u, \id_B}$ for $\nu \in
  \{ L, R \}$. We then have
  \begin{description}
    \item[{\tmem{n1}}] $\varepsilon^R_{A, v} \bullet \left( G ( \theta ) \circ
    \varepsilon_A \right) = e_{\theta} \bullet \varepsilon^R_{A, u}$;
    
    \item[{\tmem{n2}}] $\varepsilon^L_{A, v} \bullet \varepsilon_{\theta} =
    \varepsilon_{\theta} \bullet \varepsilon^L_{A, u}$;
    
    \item[{\tmem{c1}}] $\varepsilon^R_{u, B} \bullet \varepsilon^R_{A, u} =
    \varepsilon^R_{A, u}$;
    
    \item[{\tmem{c2}}] $\varepsilon^L_{A, u} \bullet \varepsilon^L_{u, B} =
    \varepsilon^L_{u, B}$;
    
    \item[{\tmem{c3}}] $\varepsilon^R_{u, B} \bullet \varepsilon^L_{A, u} =
    \varepsilon^L_{A, u} \bullet \varepsilon^R_{u, B}$.
  \end{description}
\end{theremark}

\begin{definition}
  Let $F : \as \longrightarrow \bs$ be a 2-functor. The {\tmem{identity
  distortion}} $\id_F : F \rightsquigarrow F$ is given by
  \begin{enumerateroman}
    \item $( \id_F )_f  \deq F ( f )$ for all $f \in \as_1$;
    
    \item $( \id_F )_{\alpha}  \deq F ( \alpha )$ for all $\alpha \in \as_2$;
    
    \item $( \id_F )^L_{f, g} = ( \id_F )^R_{f, g} = \id$ for all composable
    $f, g \in \as_1$.
  \end{enumerateroman}
\end{definition}

\subsection{Skew Immersions}

\begin{definition}
  An inclusion $J : \hspace{0.25em} \as \rightarrowtail \bs$ of 2-categories
  is a \textit{skew immersion} provided
  \begin{enumerate}
    \item $\as$ is a left ideal;
    
    \item there is a right ideal $\ws \subseteq \bs$ such that the
    corestriction $J : \hspace{0.25em} \as \rightarrowtail \ws$ admits a
    retraction $R_J : \hspace{0.25em} \ws \twoheadrightarrow \as$ and a
    distortion $\varepsilon : \hspace{0.25em} J \circ R_J \rightsquigarrow
    \id_{\ws}$ with $\varepsilon J = \id_{_J}$.
  \end{enumerate}
\end{definition}

\begin{theremark}
  \label{rem:sdr-skewimm}It follows by remark \ref{rem-simho} that $N_2 \left(
  \as \right)$ is a strong deformation retract of $N_2 \left( \ws \right)$
  with respect to the standard model structure on {\sset}.
\end{theremark}

For the remaining of this section, we fix a skew immersion $J : \as
\rightarrowtail \bs$ and a pushout square

\begin{center}
  $
\xy
\xymatrix{ \mathcal{A} \ar[r]^{U} 
\ar@{ >->}[d]_{J} & \mathcal{A'} \ar@{ >->}[d]^{J'} 
\\ 
\mathcal{B} \ar[r]_{W} &\mathcal{B'}} 
\POS( 8.8 , -8.8);\POS( 8.8,-11.1 ) \connect@{-}
\POS( 8.8 , -8.8);\POS( 11.1,-8.8 ) \connect@{-}
\endxy
$

\end{center}

{\noindent}along with its decomposition

\begin{center}
  $
\xy
\xymatrix{ \mathcal{A} \ar[r]^{U} \ar@{ >->}[d]_{J} 
&\mathcal{A'} \ar@{ >->}[d]^{J'} 
\\ \mathcal{W} \ar[r]_{V} \ar@{ >->}[d]_{K} & \mathcal{W'} 
\ar@{ >->}[d]^{K'} 
\\ \mathcal{B} \ar[r]_{W}&\mathcal{B'} }
\POS( 10 , -9);\POS( 10,-11.3 ) \connect@{-}
\POS( 10 , -9);\POS( 12.3,-9 ) \connect@{-}
\POS( 10 , -22);\POS( 10,-24.3 ) \connect@{-}
\POS( 10 , -22);\POS( 12.3,-22 ) \connect@{-}
\endxy
$

\end{center}

\begin{proposition}
  \label{pro:slimmersions-pstable}Skew immersions are stable under pushout.
\end{proposition}

\begin{theproof}
  By lemma \ref{lem-idealstab}, $\mathcal{A}^{\prime}$ is a left ideal and
  $\mathcal{W}^{\prime}$ is a right ideal. In particular, $\as$ is a left
  ideal in $\ws$ while $\as'$ is a left ideal in $\ws'$.
  
  Let $P$ be the 2-graph given by
\[
\xy
\xymatrix{
\forg(\as) \ar[r]^{\forg(U)} \ar@{ >->}[d]
& \forg(\as') \ar@{ >->}[d]^{\iota_{\as'}}
\\
\forg(\ws) \ar[r]_{\iota_\ws} & P
}
\POS( 13 , -9);\POS( 13,-11.4 ) \connect@{-}
\POS( 13 , -9);\POS( 15.4,-9 ) \connect@{-}
\endxy
\]
  {\noindent}The colimits in a functor category being calculated pointwise, we
  have
  \[ P \cong \left( \as'_2 + ( \ws_2 \backslash \as_2 ), \as'_1 + ( \ws_1
     \backslash \as_1 ), \as'_0 + ( \ws_0 \backslash \as_0 ) \right) \]
  with the copairs
  \[ \partial^i_P = \left[ \inn_{\as'_i} \circ \partial^i_{\as'}, \left( U_i +
     \id_{\ws_i \backslash \as_i} \right) \circ \partial^i_{\ws} |_{_{\ws_{i +
     1} \backslash \as_{i + 1}}} \right] \]
  as structural maps, for $i \in \{ 0, 1 \}$ and $\partial \in \{ \dom, \cod
  \}$. The coprojections are
  \[ \iota_{\as'} = \left( \inn_{\as'_2}, \inn_{\as'_1}, \inn_{\as'_0} \right)
  \]
  respectively
  \[ \iota_{\ws } = \left( U_2 + \id_{\bs_2 \backslash \as_2}, U_1 +
     \id_{\bs_1 \backslash \as_1}, U_0 + \id_{\bs_0 \backslash \as_0}  \right)
  \]
  Let $J' : \as' \rightarrowtail \free ( P )$and $W : \ws \rightarrow \free (
  P )$ be the morphisms of 2-graphs induced by $\iota_{\as'}$ respectively by
  $\iota_{\ws}$. Then
  \[ \ws' \cong \free ( P ) / \sim \]
  with $\sim$ the smallest sesquicongruence making $J'$ and $W$ 2-functorial
  (c.f. proposition \ref{pro:2cat-complete} and remark \ref{rem-pushinc}).
  
  Since $\as \subseteq \ws$ is a left ideal, a morphism $k$ generating $\ws$'
  has one of the following types:
    \begin{enumeratenumeric}
      \item $k \in \as'_1$;
      
      \item $k \in ( \ws_1 \backslash \as_1 ) \backslash ( \ws \backslash \as
      )_1$;
      
      \item $k \in ( \ws \backslash \as )_1$ (c.f. lemma \ref{lem:vwb});
    \end{enumeratenumeric}
  while a 2-cell $\varpi$ generating $\ws$' has one of the following types:
    \begin{enumeratenumeric}
      \item $\dom^1 ( \varpi ), \cod^1 ( \varpi ) \in \as'_1$;
      
      \item $\dom^1 ( \varpi ), \cod^1 ( \varpi ) \in ( \ws_1 \backslash \as_1
      ) \backslash ( \ws \backslash \as )_1$;
      
      \item $\dom^1 ( \varpi ), \cod^1 ( \varpi ) \in ( \ws \backslash \as
      )_1$ (c.f. lemma \ref{lem:vwb}).
    \end{enumeratenumeric}
  In particular, given a morphism $k$ of type 2 we have $\dom ( k ) = U ( A )$
  for some $A \in \as_0$ and $\cod ( k ) = B$ for some $B \in \ws_0 \backslash
  \as_0$. A typical situation can be depicted as follows

  \begin{center}
    $
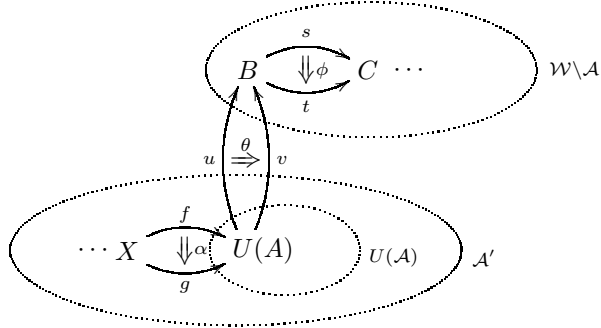
\xy
\xymatrix @*[r] {
 && B \rtwocell^s_t{\phi} & C 
\\
&&&
\\
&X \rtwocell^f_g{\alpha} & U(A) \uutwocell^u_v{\theta} &
}
\POS(33.5,-24)
\drop{\xycircle(10,6){.}}
\POS(48,-25)
\drop{\scriptstyle U(\as)}
\POS(27,-24)
\drop{\xycircle(30,10){.}}
\POS(60,-25)
\drop{\scriptstyle \as'}
\POS(45,0)
\drop{\xycircle(22,9){.}}
\POS(72,0)
\drop{\scriptstyle \ws\setminus\as}
\POS(8,-24)
\drop{\cdots}
\POS(50,0)
\drop{\cdots}
\endxy
$

  \end{center}

  {\noindent}General morphisms of $\ws'$ are thus composable strings
  \[ < f_0 ; \cdots ; f_n ; u ; s_1 ; \cdots ; s_m > \]
  with $f_0, \ldots, f_n$ of type 1, $u$ of type 2 and $s_1, \ldots, s_m$ of
  type 3. Simlarly, general 2-cells of $\ws'$ are horizontally composable
  strings
  \[ \ll \alpha_0 ; \cdots ; \alpha_n ; \theta ; \phi_1 ; \cdots ; \phi_m \gg
  \]
  with $\alpha_0, \ldots, \alpha_n$ of type 1, $\theta$ of type 2 and $\phi_1,
  \ldots, \phi_m$ of type 3. On the other hand, the relations governing $\ws'$
  impose the identities
  \[ < f_0 ; \cdots ; f_n ; u ; s_1 ; \cdots ; s_m > = < f_n \circ \cdots
     \circ f_0 > ; < s_m \circ \cdots \circ s_1 \circ u > \]
  respectively
  \[ \ll \alpha_0 ; \cdots ; \alpha_n ; \theta ; \phi_1 ; \cdots ; \phi_m \gg
     = \ll \alpha_n \circ \cdots \circ \alpha_0 \gg ; \ll \phi_m \circ \cdots
     \circ \phi_1 \circ \theta \gg \]
  among such strings.
  
  Finally, there is the retraction
  \[ R'_{J'} = \left( \left[ \id_{\as'_2}, ( U \circ R )_2 \right], \left[
     \id_{\as'_1}, ( U \circ R )_1 \right], \left[ \id_{\as'_0}, ( U \circ R
     )_0 \right] \right) : \ws' \longrightarrow \as' \]
  given by universal property. Let $\nu \in \{ L, R \}$. There is the
  distortion
  \[ \xi : J' \circ R'_{J'} \rightsquigarrow \id_{\ws'} \]
  given by
  \begin{itemizeminus}
    \item $\xi_{< f >}  \deq < f >$, $\xi_{\ll \alpha \gg}  \deq \ll \alpha
    \gg$  and $\xi^{\nu}_{< f >, < g >}  \deq \ll \gg_{< g \circ f >}$ for all
    $f$ of type 1, all $\alpha$ of type 1 respectively all composable $f$ and
    $g$ of type 1;
    
    \item $\xi_{< p >}  \deq < \varepsilon_p$>, $\xi_{\ll \beta \gg}  \deq \ll
    \varepsilon_{\beta} \gg$ and $\xi^{\nu}_{< p >, < q >}  \deq \ll
    \varepsilon^{\nu}_{p, q} \gg$for all $p$ of type 2 and 3, all $\beta$ of
    type 2 and 3 respectively all composable $p$ and $q$ of type 2 or 3;
    
    \item $\xi_{< f ; u >}  \deq < f ; \varepsilon_u >$ for all composable $f$
    of type 1 and $u$ of type 2;
    
    \item $\xi_{\ll \alpha ; \theta \gg}  \deq \ll \alpha ;
    \varepsilon_{\theta} \gg$for all horizontally composable $\alpha$ of type
    1 and $\theta$ of type 2;
    
    \item $\xi^{\nu}_{\dom ( f ), < f, u >}  \deq \ll f ; \xi^{\nu}_{U ( \dom
    ( f ) ), u} \gg = \ll f ; \varepsilon^{\nu}_{\dom ( f ), u} \gg$ and
    
    \item $\xi^{\nu}_{< f, u >, \cod ( u )}  \deq \ll f ; \xi^{\nu}_{u, \cod (
    u )} \gg = \ll f ; \varepsilon^{\nu}_{u, \cod ( u )} \gg$ for all
    composable $f$ of type 1 and $u$ of type 2.
  \end{itemizeminus}
  The axioms of distortion are easily checked, e.g.
  \begin{eqnarray*}
    \ll \xi_X ; \ll \alpha ; \theta \gg \gg : \xi^R_{X, < g ; v >} & = & \ll
    \alpha ; \theta \gg : \xi^R_{X, < g ; v >} ( \xi_X = \id_X )\\
    & = & \ll \alpha ; \theta \gg : \ll g ; \xi^R_{U ( A ), v} \gg\\
    & = & \ll \alpha \gg ; \ll \theta : \xi^R_{U ( A ), v} \gg\\
    & = & \ll \alpha ; \varepsilon^R_{A, v} \bullet \theta \gg\\
    & = & \ll \alpha ; \varepsilon^R_{A, v} \bullet ( \theta \circ
    \varepsilon_A ) \gg ( \varepsilon_A = \id_A )\\
    & = & \ll \alpha ; \varepsilon_{\theta} \bullet \varepsilon^R_{A, u}
    \gg\\
    & = & \ll \alpha \gg ; \ll \xi^R_{U ( A ), u} : \varepsilon_{\theta}
    \gg\\
    & = & \ll f ; \xi^R_{U ( A ), u} \gg : \ll \alpha ; \varepsilon_{\theta}
    \gg\\
    & = & \xi^R_{X, < f, u >} : \xi_{\ll \alpha, \theta \gg}
  \end{eqnarray*}
  (c.f. remark \ref{rem-dist-conv}), while $\xi J' = \id_{J'}$ holds by
  construction.
\end{theproof}

\subsection{The 2-Thomason Model Structure}

\begin{lemma}
  \label{lem:vvprime} $\left( \bs \backslash \as \right) \cap \ws \cong \left(
  \bs' \backslash \as' \right) \cap \ws'$.
\end{lemma}

\begin{theproof}
  Let $\star_{\as} \in \bs / \as$ be the object such that $( \rho_{\as} \circ
  !_{\as} ) ( \as ) = \left\lceil \ast_{\as} \right\rceil$ (c.f. defintion
  \ref{def-thoquot}). Given the iso $i : \bs / \as \cong \bs' / \as'$ (c.f.
  proposition \ref{prop-thoquot}), it is immediate that $i ( \star_{\as} ) =
  \star_{\as'}$. On the other hand, $\cod ( f ) \in \ws_1$ for all $\left. f
  \in \left( \bs_1 \backslash \as_1 \right) \backslash ( \bs \backslash \as
  \right)_1$ since $\ws$ is a right ideal and $\cod ( f' ) \in \ws'_1$ for all
  $\left. f' \in \left( \bs'_1 \backslash \as'_1 \right) \backslash ( \bs'
  \backslash \as' \right)_1$ since $\ws'$ is a right ideal. Hence there is a
  bijection
  \[ \left( \left( \bs \backslash \as \right) \cap \ws \right)_0 \cong \left(
     \left( \bs' \backslash \as' \right) \cap \ws' \right)_0 \]
  induced by $i$. But $\left( \bs \backslash \as \right) \cap \ws$ is a right
  ideal in $\bs \backslash \as$ and $\left( \bs' \backslash \as' \right) \cap
  \ws'$ is a right ideal in $\bs' \backslash \as'$. In particular, both
  sub-2-categories are full and locally full, hence $i_{\text{$\left( \bs
  \backslash \as \right) \cap \ws$}}$ is an iso. 
\end{theproof}

\begin{definition}
  Let $\mathbb{M}$ be a model category. A weak pushout square in $\mathbb{M}$
  is a commuting square such that the comparison map from the inscribed
  pushout is a weak equivalence:

  \begin{center}
    $
\xy
\xymatrix{
A \ar[d] \ar[r] 
&
B \ar[d] \ar@/^1pc/[ddr]
&
\\
C \ar[r]  \ar@/_1pc/[drr] & 
B +_A C \ar@{.>}[dr]^{\sim}
&
\\
&& D
}
\POS( 13 , -8);\POS( 13,-10.3 ) \connect@{-}
\POS( 13 , -8);\POS( 15.3,-8 ) \connect@{-}
\endxy
$

  \end{center}
\end{definition}

\begin{lemma}
  \label{lem:homopush}The image under $N_2$ of a pushout square of a skew
  immersion along an arbitrary 2-functor is a weak pushout square.
\end{lemma}

\begin{theproof}
  Consider
  
  \begin{center}
    $
 \xymatrix @*[l]{
   & 
   \ar@{}[dr]|{(1)}
   \mathcal{A} \ar[r]^{U} \ar@{ >->}[d]_{J} & 
   \mathcal{A'} \ar@{ >->}[d]^{J'} 
 \\ 
   \ar@{}[dr]|{(2)}
(\bs\setminus\as) \cap \mathcal{W} 
\ar@{ >->}[r]^{\omega_\ws} \ar@{ >->}[d] & 
   \ar@{}[dr]|{(3)}
   \mathcal{W} \ar[r]^{V} \ar@{ >->}[d]_{K} & 
   \mathcal{W'} \ar@{ >->}[d]^{K'} 
 \\ 
\bs\setminus\as \ar@{ >->}[r]_{\omega_\bs} &
   \mathcal{B} \ar[r]_{W} & 
   \mathcal{B'}
 }
$

  \end{center}

  {\noindent}By remark \ref{rem:sdr-skewimm}, $N_2 \left( J \right)$ is part
  of a deformation-retraction in $\sset$ and hence an acyclic cofibration.
  $N_2 \left( J' \right)$ is an acyclic cofibration for the same reason. Since
  the latter are preserved by pushouts in any model category, $N_2$ carries
  square $\left( 1 \right)$ to a weak pushout square by 2-of-3.
  
  On the other hand, $N_2$ carries square $\left( 2 \right)$ to a pushout
  square by lemma \ref{lem:vwb}. Now $\bs \backslash \as \cong \bs' \backslash
  \as'$ by corollary \ref{cor-thoquot} and $\left( \bs \backslash \as \right)
  \cap \ws \cong \left( \bs' \backslash \as' \right) \cap \ws'$ by lemma
  \ref{lem:vvprime}. Moreover,
  \[ V \circ \omega_{\ws} = \omega_{\ws'} \]
  and
  \[ W \circ \omega_{\bs} = \omega_{\bs'} \]
  by construction of the pushout squares (1) and (3) (c.f. remark
  \ref{rem-pushinc}). Hence the joint square $\left( 2 \right) \left( 3
  \right)$ also becomes a pushout square under $N_2$ by lemma \ref{lem:vwb}.
  
  But then square $\left( 3 \right)$ is also transformed in a pushout square
  by $N_2$ as a consequence of the pushout lemma and the assertion follows
  applying the glueing lemma.
\end{theproof}

\begin{lemma}
  \label{lem:sievecond}Let $A \subseteq B$ be an inclusion of posets. $C_2 N_1
  \left( A \hookrightarrow B \right)$ is a left ideal if $A$ is downclosed and
  is a right ideal if $A$ is upper-closed.
\end{lemma}

\begin{theproof}
  Immediate.
\end{theproof}

\begin{lemma}
  \label{lem:wimcond}Let $A \subseteq B$ be an inclusion of posets with $A$
  down-closed. Let $\upa$ be $A$'s upper-closure. If $i : \hspace{0.25em} A
  \hookrightarrow \hspace{0.25em} \upa$ admits a retraction $r$ such that $( i
  \circ r ) ( x ) \leq x$ for all $x \in \hspace{0.25em} \upa$, then $C_2 N_1
  \left( A \hookrightarrow B \right)$ is a skew immersion.
\end{lemma}

\begin{theproof}
  By lemma \ref{lem:sievecond}, $C_2 N_1 \left( A \right)$ is a left ideal and
  $C_2 N_1 \left( \upa \right)$ a right ideal. It is easy to see that the
  inclusion $\overline{i} : C_2 N_1 \left( A \right) \hookrightarrow C_2 N_1
  \left( \upa \right)$ admits a retraction $\overline{r} : \hspace{0.25em} C_2
  N_1 \left( \upa \right) \rightarrow C_2 N_1 \left( A \right)$ and that there
  is the family $\left( \varepsilon_x : \hspace{0.25em} \left( \overline{i}
  \circ \overline{r} \right) x \rightarrow x \right)$ given by the
  inequalities $( i \circ r ) ( x ) \leq x$. Since the 2-categorification of a
  poset is a locally ordered 2-category, this family determines a distortion. 
\end{theproof}

\begin{lemma}
  \label{lem:sdhorncalc}Let $f : \hspace{0.25em} \ord \rightarrow \ord$ be the
  functor assigning to an order the order of its non-empty totally ordered
  finite subsets, ordered by inclusion. Let
  \[ H_{k, n} \deq f \left( \left[ n \right] \right) \setminus \left\{ \left(
     0, \ldots n \right), \left( 0, \ldots, k - 1, k + 1, \ldots, n \right)
     \right\} \]
  Then
  \[ Sd^2 \left( \Lambda^k \left[ n \right] \right) = N_1 \left( f \left(
     H_{k, n} \right) \right) \]
  and
  \[ Sd^2 \left( \Delta \left[ n \right] \right) = N_1 \left( f^2 \left(
     \left[ n \right] \right) \right) \]
\end{lemma}

\begin{theproof} The subdivision of a simplicial complex is the nerve of its
poset of non-degenerate faces.
\end{theproof}

\begin{lemma}
  \label{lem:collar}Let $f$ be as in lemma \ref{lem:sdhorncalc} and $P$ be a
  finite connected poset with a greatest element $\top$. Let further $k \in P$
  be a maximal element of $P \setminus \top$ and $P_k \deq P \setminus \left\{
  k, \top \right\}$. Finally, let $P$'s $k$-\textit{horn} be given by $H_{P,
  k} \deq f \left( P_k \right)$ and $P$'s $k$-\textit{collar} be given by
  $C_{P, k} \deq f \left( P \right) \setminus \left\{ \left( \top \right),
  \left( k \right), \left( k, \top \right) \right\}$. Then
  \[ \ups{H_{P, k}} = C_{P, k} \]
  and the assignment
  \[ \begin{array}{rlcl}
       r : & C_{P, k} & \rightarrow & H_{P, k}\\
       & x & \mapsto & \max \left( \left( \downarrow x \right) \cap H_{P, k}
       \right)
     \end{array} \]
  determines a retraction such that $r ( x ) \subseteq x$ for all $x \in C_{P,
  k}$.
\end{lemma}

\begin{theproof}
  By the very definition, for any list $x \in C_{P, k}$ there is a list
  $x^{\prime} \in H_{P, k}$ such that $x^{\prime} \subseteq x$, hence the
  first assertion and also $\emptyset \neq \left( \downarrow x \right) \cap
  H_{P, k}$. The latter has a greatest element for
  \begin{enumerate}
    \item if $x \in H_{P, k}$ then $\max \left( \left( \downarrow x \right)
    \cap H_{P, k} \right) = x$;
    
    \item if $x \nin H_{P, k}$ then, by hypothesis on $k$, there are lists
    $x^{\prime} \in H_{P, k}$ and $x^{\prime \prime} \in \left\{ \left( \top
    \right) \left( k \right), \left( k, \top \right) \right\}$ such that $x$
    is the concatenation $x = x^{\prime} \ast x^{\prime \prime}$ hence $\max
    \left( \left( \downarrow x \right) \cap H_{P, k} \right) = x^{\prime}$.
  \end{enumerate}
\end{theproof}

\begin{lemma}
  \label{lem-hornskim}Let $i_{k, n} : \hspace{0.25em} \Lambda^k \left[ n
  \right] \rightarrow \Delta \left[ n \right]$ be a horn inclusion. Then $C_2
  \left( \tmop{Sd}^2 \left( i_{k, n} \right) \right)$ is a skew immersion.
\end{lemma}

\begin{theproof}
  Clearly, $f \left( H_{k, n} \right) \subseteq f^2 \left( \left[ n \right]
  \right)$ is down-closed. The assertion readily follows by lemma
  \ref{lem:sdhorncalc}, lemma \ref{lem:collar} and lemma \ref{lem:wimcond}. 
\end{theproof}

\begin{theorem}
  \label{thm:model} $Ex^2 \circ N_2$ creates a model structure on 2-Cat.
\end{theorem}

\begin{theproof}
  We need to show that the conditions of proposition \ref{prop-creation}  
are satisfied.
  Condition (i) holds since {\tmstrong{{\tcat}}} is locally finitely
  presentable (c.f. remark \ref{rem-2catlocpres}). Condition (ii) is a
  well-known fact about the standard model category structure on $\sset$.
  
  To verify condition (iii), observe that for any ordinal $\lambda$ and any
  $\lambda$-sequence $X : \lambda \longrightarrow \text{{\tmstrong{{\tcat}}}}$
  \begin{eqnarray*}
    ( \ex^2 \circ N_2 ) ( \colim_{\lambda} X )_n & = \sset \left( \Delta [ n
    ], ( \ex^2 \circ N_2 ) ( \colim_{\lambda} X ) \right)\\
    & \cong \sset \left( ( C_2 \circ \sd^2 ) ( \Delta [ n ] ),
    \colim_{\lambda} X \right)\\
    & \cong \colim_{\lambda}  \sset \left( ( C_2 \circ \sd^2 ) ( \Delta [ n ]
    ), X ) \right.\\
    & \cong \colim_{\lambda}  \sset \left( \Delta [ n ], ( \ex^2 \circ N_2 )
    ( X ) \right)\\
    & \cong \colim_{\lambda} ( \ex^2 \circ N_2 ) ( X )_n
  \end{eqnarray*}
  for all $n \geq 0$. The third equality is due to the fact that $\sset$, as
  any topos of presheaves, is ($\aleph_0$-) locally presentable so in
  particular simplicial sets are small with respect to the class of all
  simplicial morphisms. Since colimits are calculated dimension-wise in
  $\sset$, it follows that $\ex^2 \circ N_2$ commutes with colimits of
  $\lambda$-sequences.
  
  To complete the proof we must show that for any pushout diagram in
  {\tmstrong{{\tcat}}}
\[
\xy
\xymatrix{
(C_2 \circ \sd^2)(\Lambda^k[n]) \ar[d]_{(C_2 \circ \sd^2) (j_{n,k})} 
\ar [r]^(0.7)f
&\as \ar[d]^g
\\
(C_2 \circ \sd^2)(\Delta[n]) \ar[r]_(0.7){\bar f} & \bs
}

\POS( 18.5 , -10);\POS( 18.5,-12.4 ) \connect@{-}
\POS( 18.5 , -10);\POS( 20.9,-10 ) \connect@{-}
\endxy
\]
  {\noindent}for any $n > 0$ and $0 \leq k \leq n$, $( \ex^2 \circ N_2 ) ( g
  )$ is a weak equivalence of simplicial sets. Consider
\[
\xy
\xymatrix{
\sd^2 (\Lambda ^k[n]) \ar[r]^(.4){\eta_{\Lambda ^k[n]}}
\ar[d]_{\sd^2 (j_{n,k})}
&
(N_2\circ C_2 \circ \sd^2)( \Lambda ^k[n])
\ar [r]^(.7){ N_2(f)}
\ar [d]_{ (N_2 \circ C_2 \circ \sd^2) (j_{n,k})} 
&
N_2(\as) \ar[d]^\varphi \ar@/^3ex/[dr]^{ N_2(g)}
&
\\
\sd^2(\Delta [n]) \ar [r]_(.4){\eta_{\Delta[n]}}
&
(N_2 \circ C_2 \circ \sd^2)( \Delta[n])
\ar [r]_(.7)\psi 
\ar@/_6ex/[rr]_{N_2(\bar{f})}
&
\bullet \ar@{.>}[r]^\omega
& 
N_2(\bs)
}
\POS( 62 , -10);\POS( 62,-12.4 ) \connect@{-}
\POS( 62 , -10);\POS( 64.4,-10 ) \connect@{-}
\endxy
\]
  {\noindent}with  $\omega$ the comparison morphism. Since $\sd^2 ( K )$ is
  the $1$-nerve of a poset for any simplicial set $K$, the unit maps $\eta$
  are isos by remark \ref{rem-n2c2n1}, so in particular weak equivalences.
  Furthermore, there is the obvious commutative diagram
  \begin{center}
    $
\xymatrix{ 
\Lambda ^k[n]\ar [r]\ar [d]_{j_{n,k}}
&
\sd^2 \Lambda ^k[n]\ar [d]^{\sd^2 (j_{n,k})}
\\
\Delta [n]\ar [r]
&
\sd^2\Delta [n]
}$

  \end{center}
  {\noindent}in which the horizontal maps induce homeomorphisms after
  geometric realization and are therefore weak equivalences. Hence, by 2-of-3,
  $\sd^2 ( j_{n, k} )$ is also a weak equivalence. Thus, applying 2-of-3 to
  the lefthand square of diagram (*), we obtain that $( N_2 \circ C_2 \circ
  \sd^2 ) ( j_{n, k} )$ is a weak equivalence as well, which implies that
  $\varphi$ is a weak equivalence, since acyclic cofibrations are preserved
  under pushout in any model category.
  
  On the other hand, $\omega$ is also a weak equivalence, as $( C_2 \circ
  \sd^2 ) ( j_{n, k} )$ is a skew immersion by lemma \ref{lem-hornskim}. Thus,
  $N_2 ( g ) = \omega \circ \varphi$ is a weak equivalence, which implies that
  $( \ex^2 \circ N_2 ) ( g )$ is a weak equivalence since $\ex$ preserves the
  latter, which completes the proof.
\end{theproof}

We call the model structure of theorem \ref{thm:model} the \textit{2-Thomason
model structure} since it is conceptually similar to the model structure on
\textbf{$\cat$} due to R.W.Thomason {\cite{thoma-model}}.

\section{\label{sec:Homotopy}Homotopy}

\begin{definition}
  Let $\as$ be a 2-category and $f, g \in \as_1$.
  \begin{enumerate}
    \item A lax square $\left( u_0, u_1, \alpha \right) : f \longrightarrow g$
    is given by the diagram

    \begin{center}
      $\xymatrix{
X \ar[r]^{u_0} \ar[d]_{f} 
\drtwocell<\omit>{^{\alpha}}
& Y \ar[d]^{g} 
\\
X' \ar[r]_{u_1} & 
Y'
}$

    \end{center}
    
    {\noindent}Let $( v_0, v_1, \beta ) : g \longrightarrow h$ be a further
    lax square. Their {\tmem{pasting composite}} is the lax square
    \[ \left. ( v_0, v_1, \beta ) \circledast ( u_0, u_1, \alpha \right)  \deq
       \left ( v_0 \circ u_0, v_1 \circ u_1, ( \beta \circ u_0 ) \bullet ( v_1 \circ
       \alpha ) \right ): f \longrightarrow h \]
    \item A cylinder $\left( \theta_0, \theta_1 \right) : \left( u_0, u_1,
    \alpha \right) \longrightarrow \left( v_0, v_1, \beta \right)$ is given by
    the diagram

    \begin{center}
      $
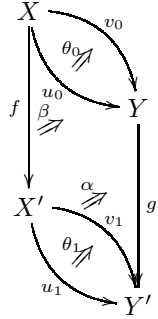
\xymatrix{ 
X \ar[dd]_{f} 
\drtwocell<5>^{v_0}_{u_0}{^{\theta_0}} 
\ddrtwocell<\omit>{^{<2.8>\beta}} & 
\\
\ddrtwocell<\omit>{^{<-1>\alpha}} & Y \ar[dd]^{g} 
\\
X' \drtwocell<5>^{v_1}_{u_1}{^{\theta_1}} & 
\\
& Y'
}$

    \end{center}

    {\noindent}where $\left( g \hc \theta_0 \right) \bullet \alpha = \beta
    \bullet \left( \theta_1 \hc f \right)$
  \end{enumerate}
\end{definition}

\begin{proposition}
  \label{pro:cyl}Let $\as$ be a 2-category.
  \begin{enumeratenumeric}
    \item There is a 2-category $\cyl( \as )$ \textit{}given by the data
    \begin{itemizeminus}
      \item Objects: morphisms of $\as_{}$;
      
      \item Morphisms: lax squares;
      
      \item 2-cells: cylinders.
    \end{itemizeminus}
    Composition of morphisms is given by pasting while the operations on
    2-cells are those of $\as$ taken componentwise.
    
    \item The assignments
    \[ \begin{array}{llll}
         \dom_{\as} : & \cyl ( \as ) & \longrightarrow & \as\\
         & f & \longmapsto & \dom ( f )\\
         & \left( u_0, u_1, \alpha \right) & \longmapsto & u_0\\
         & ( \theta_0, \theta_1 ) & \longmapsto & \theta_0
       \end{array} \]
    \[ \begin{array}{llll}
         \cod_{\as} : & \cyl ( \as ) & \longrightarrow & \as\\
         & f & \longmapsto & \cod ( f )\\
         & \left( u_0, u_1, \alpha \right) & \longmapsto & u_1\\
         & ( \theta_0, \theta_1 ) & \longmapsto & \theta_1
       \end{array} \]
    and
    \[ \begin{array}{llll}
         I_{\as} : & \as & \longrightarrow & \cyl \left( \as \right)\\
         & X & \longmapsto & \id_X\\
         & f & \longmapsto & ( f, f, \id_f )\\
         & \alpha & \longmapsto & ( \alpha, \alpha )
       \end{array} \]
    are 2-functorial.
  \end{enumeratenumeric}
  
\end{proposition}

Following B\'enabou, we 
call $\cyl \left( \as \right)$ the \textit{2-category of
cylinders over} $\as$ \cite{bena}. The name stems from the
``geometry'' of 2-cells. Notice that $\cyl \left( \as \right)$ is a
generalization of the familiar category of arrows. The construction is
2-functorial, yet this fact is not relevant for the present development.

\begin{theremark}
  \label{rem-nlaxtocyl}Let $\as$ be a 2-category and $q \in \N$. A normal lax
  functor $F : [ q ] \longrightarrow \cyl ( \as )$ is determined by
  \begin{itemizeminus}
    \item a morphism $F ( k ) : F ( k )^- \longrightarrow F ( k )^+$ for all
    $0 \leqslant k \leqslant q$;
    
    \item a lax square $F ( k < l ) \deq  \left( F ( k, l )^-, F ( k, l )^+, F
    ( k, l ) \right) : F ( k ) \longrightarrow F ( l )$ as in

    \begin{flushleft}
      \begin{center}
        $\xymatrix{
F(k)^- \ar[rr]^{F(k,l)^-} \ar[d]_{F(k)} 
\POS( 10, -10);\POS( 17,-6) \connect@{=>}
\POS(12,-4.5) \drop{\scriptstyle F(k,l)}
& & F(l)^- \ar[d]^{F(l)} 
\\
F(k)^+ \ar[rr]_{F(k,l)^+} & &
F(l)^+
}$

      \end{center}
    \end{flushleft}
    
    {\noindent}for all $0 \leqslant k < l \leqslant q$;
    
    \item a cylinder
    \[ F ( k < l < m ) \deq  \left( ( F ( k, l, m )^-, F ( k, l, m )^+ \right)
       : F ( l < m ) \circledast F ( k < l ) \longrightarrow F ( k < m )^{} \]
    for all $0 \leqslant k < l < m \leqslant q$
  \end{itemizeminus}
  such that
  \[ F ( k, m, n )^s \bullet \left( F ( m, n )^s \circ F ( k, l, m )^s \right)
     = F ( k, l, n )^s \bullet \left( F ( l, m, n )^s \circ F ( k, l )^s
     \right) \]
  for all $s \in \{ -, + \}$ and $0 \leqslant k < l < m < n \leqslant q$.
\end{theremark}

\begin{definition}
  Let $F, G : \hspace{0.25em} \as \rightarrow \mathcal{B}$ be 2-functors. A
  lax transformation $\alpha : \hspace{0.25em} F \Rightarrow G$ is given by
  \begin{itemizeminus}
    \item a morphism $\mathcal{} \alpha_X : \hspace{0.25em} F \left( X \right)
    \rightarrow G \left( X \right)$ for each $X \in \as$ and
    
    \item a 2-cell

    \begin{center}
      $\xymatrix{
F(X) \ar[r]^{\alpha_X} \ar[d]_{F(f)} 
\drtwocell<\omit>{ \alpha_f}
& G(X) \ar[d]^{G(f)} 
\\
F(Y) \ar[r]_{\alpha_Y} & 
G(Y)
}$

    \end{center}
    
    {\noindent}for each morphism $f : \hspace{0.25em} X \rightarrow Y$
  \end{itemizeminus}
  such that
  \begin{enumerateroman}
    \item $\alpha_{f'} \bullet \left( G \left( \theta \right) \circ \alpha_X
    \right) = \left( \alpha_Y \circ F \left( \theta \right) \right) \bullet
    \alpha_f$  for each 2-cell $\theta : \hspace{0.25em} f \Rightarrow f' : X
    \rightarrow Y$;
    
    \item $\left( \alpha_g \circ F \left( f \right) \right) \bullet \left( G
    \left( g \right) \circ \alpha_f \right) = \alpha_{g \circ f}$  for each $f
    : \hspace{0.25em} X \rightarrow Y$ and $g : \hspace{0.25em} Y \rightarrow
    Z$.
  \end{enumerateroman}
  
\end{definition}

\begin{proposition}
  \label{pro:benabou-oplax}The following are equivalent
  \begin{enumerateroman}
    \item There is a lax transformation $\alpha : \hspace{0.25em} F
    \Rightarrow G$;
    
    \item There is a 2-functor $\overline{\alpha} : \hspace{0.25em} \as
    \rightarrow \cyl \left( \mathcal{B} \right)$ such that

    \begin{center}
      $\xymatrix{
&& \cyl ( \mathcal{B} ) 
\ar[d]^{\langle\dom_\bs,\cod_\bs\rangle} \\
\mathcal{A} \ar[urr]^{\bar{\alpha}} 
\ar[rr]_{\scriptscriptstyle \langle F,G \rangle} & &
\mathcal{B} \times \mathcal{B} 
}$

    \end{center}

    {\noindent}commutes.
  \end{enumerateroman}
  
\end{proposition}

Our 2-category of cylinders is in fact the strict case of B\'enabou's
\textit{bicategory of cylinders}. He defined lax transformations for lax
functors among bicategories in terms of this classifing device \cite{bena}.

\begin{definition}
  Let $\mathbb{M}$ be a model category and $P, B \in \M$. $P$ is a
  \textit{path object} on $B$ if there is a a morphism $p_B : P
  \longrightarrow B \times B$ and commuting diagram

  \begin{center}
    $\xymatrix{ 
&& P \ar[d]^{p_B}  \\
B \ar[urr]^{\sim} \ar[rr]_{\Delta} & & B \times B
}$

  \end{center}

\end{definition}

\begin{proposition}
  \label{pro:cyl} $\cyl \left( \mathcal{A} \right)$ is a path object on
  $\mathcal{A}$ in the 2-Thomason model structure.
\end{proposition}

\begin{theproof}
  It is immediate that
  \begin{center}
    $\xymatrix{ 
&& \cyl(\as) \ar[d]^{\langle\dom_\as,\cod_\as\rangle}  \\
\as \ar[urr]^{I_\as} \ar[rr]_{\Delta} & & \as \times \as
}$

  \end{center}
  {\noindent}commutes.
  On the other hand, there is a simplicial homotopy
  \begin{center}
    $\xymatrix{ 
N_2\left ( \cyl (\as) \right ) 
\ar[r]^(.4){i_0}
\ar[ddr]_\id
&
N_2\left ( \cyl (\as) \right ) \times [1]
\ar[dd]^H
&
N_2\left ( \cyl (\as) \right )
\ar[l]_(.4){i_1}
\ar[ddl]^{N_2 (I_\as \circ \dom_\as)}
\\
\\
& N_2\left ( \cyl (\as) \right ) &
}$

  \end{center}
  {\noindent}It can be constructed as a family
  \[ H^n_i : N_2 ( \cyl ( \as ) )_n \longrightarrow N_2 ( \cyl ( \as ) )_{n +
     1}, 0 \leqslant i \leqslant n, n \geqslant 0 \]
  enjoying the well-known properties. Let
  \[ \cart_F ( s, t ) \deq  \left( F ( s, t )^-, F ( t ) \circ F ( s, t )^-
     \right) : \id_{F ( s )^-} \longrightarrow F ( t ) \]
  and
  \[ \dom_F ( s, t ) \deq  \left( F ( s, t )^-, F ( s, t )^- \right) : \id_{F
     ( s )^-} \longrightarrow \id_{F ( t )^-} \]
  be lax squares for $0 \leqslant s < t \leqslant n$. Let $F \in N_2 ( \cyl (
  \as ) )_n \cong \tmmathbf{\nlax} \left( [ n ], \cyl ( \as ) \right)$. The
  normal lax functor
  \[ H^n_i ( F ) : [ n + 1 ] \longrightarrow \cyl ( \as ) \]
  is given by the following data:
  \begin{itemizeminus}
    \item $H^n_i ( F ) ( p ) \deq  \left\{  \begin{array}{ll}
      \id_{F ( p )^-} & p \leqslant i\\
      F ( p - 1 ) & p > i
    \end{array}\right.$
    
    \item $H^n_i ( F ) ( p < q ) \deq  \left\{  \begin{array}{ll}
      \dom_F ( p, q ) & p, q \leqslant i\\
      \cart_F ( p, q - 1 ) & p \leqslant i, q > i\\
      F ( p - 1 < q - 1 ) & p, q > i
    \end{array}\right.$
    
    \item 
$
H^n_i ( F ) ( p < q < r ) \deq  
\\ \\
\;\;\;\;\;\;\;\;
 \left\{ \begin{array}{ll}
      \left( F ( p, q, r )^-, F ( p, q, r )^- \right) & p, q, r \leqslant i\\
      \left( F ( p, q, r - 1 )^-, F ( r - 1 ) \circ F ( p, q, r - 1 )^-
      \right) & p, q \leqslant i \wedge r > i\\
      \left( F ( p, q - 1, r - 1 )^- \right. , 
 
& p \leqslant i \wedge q, r > i\\
\left( F ( r - 1 ) \circ F ( p, q - 1, r- 1 )^- \right)  \bullet 

\left. \left( F ( q - 1, r - 1 ) \circ F ( p, q - 1 )^-\right) \right)
\\
      F ( p - 1 < q - 1 < r - 1 ) & p, q, r \geqslant i
    \end{array} \right.
$
  \end{itemizeminus}
(c.f. remark \ref{rem-nlaxtocyl}). A laborious yet straightforward calculation shows that the coherence
  conditions hold and that the $H^n_i$'s commute with faces and 
degeneracies as 
required.
  It thus follows (by functoriality) that there is a homotopy
  \[ | N_2 \left( I_{\as} \circ \dom_{\as} \right) | \sim \id_{| N_2 \left(
     \cyl \left( \as \right) \right) |} \]
  hence $I_{\as}$ is a homotopy equivalence so in particular a weak
  equivalence.
\end{theproof}

\begin{definition}
  Let $\mathbb{M}$ be a model category. \textit{Given} $f, g : \hspace{0.25em}
  A \rightarrow B$, there is a \textit{right homotopy $f \simeq g$} if there
  is a path object over $B$ such that $\left\langle f, g \right\rangle$
  factors through $p_B$ as in

  \begin{center}
    $\xymatrix{ 
&& P \ar[d]^{p_B} \\
A \ar[urr] \ar[rr]_{\langle f,g \rangle} & & B \times B
}$

  \end{center}

\end{definition}

\begin{theorem}
  Lax transformations are right homotopies in the 2-Thomason model structure.
\end{theorem}

\begin{theproof}
  Direct consequence of proposition \ref{pro:cyl}.
\end{theproof}

\begin{theremark}
  Reversing the direction of the 2-cell in the definition of a lax square
  yields the dual notion of {\tmem{oplax square}} and those of
  {\tmem{opcylinder}} and of {\tmem{oplax transformation}} respectively. It is
  easy to see that oplax cylinders  are path objects and, consequently, oplax
  transformations are right homotopies in the 2-Thomason model structure.
\end{theremark}

\bibliographystyle{abbrv}
\bibliography{ProcessSemantics,beke,benabou,cat,catta,cisinski,gab,god,goub,gray,hermida,kw,milner,modelcats,street,winskel,wolff}

\end{document}